\documentclass[]{compositio}

\usepackage{amsbsy,amssymb,amscd,amsfonts,latexsym,amstext,delarray,amsmath,stackrel,lineno,tabularx,url,MnSymbol,tikz,cite,yfonts}
\usepackage{float} 
\usepackage[pdfborder={0 0 0}]{hyperref}
\usepackage{mathrsfs}
\usepackage[all]{xy}
\hypersetup{
    colorlinks = true,
    linkcolor = blue,
    anchorcolor = red,
    citecolor = blue,
    filecolor = red,
    pagecolor = red,
    urlcolor = magenta
}

\usepackage{graphicx}

\usepackage{pdfpages}

\usepackage{blindtext}
\usepackage{etoolbox}
\newcommand{\printthis}[2][true]{%
\ifbool{#1}{%
#2%
}{%
}%
}

\newtheorem{thm}{Theorem}[section] 
\newtheorem{defn}[thm]{Definition} 
\newtheorem{prop}[thm]{Proposition}
\newtheorem{cor}[thm]{Corollary}
\newtheorem{lem}[thm]{Lemma}

\newtheorem{rem}[thm]{Remark}

\def\Hom{{\rm Hom}}

\def\Spec{{\rm Spec\,}}
\def\root{{\rm Root\,}}

\def\Mo{\mathfrak{Mo}}
\def\A{{\mathbb A}}

\def\F{{\mathbb F}}

\def\N{{\mathbb N}}

\def\Q{{\mathbb Q}}
\def\R{{\mathbb R}}
\def\Z{{\mathbb Z}}

\def\fC{{\mathfrak C}}

\def\cB{{\mathcal B}}
\def\cC{{\mathcal C}}
\def\cD{{\mathcal D}}
\def\cE{{\mathcal E}}
\def\cG{{\mathcal G}}

\def\cF{{\mathcal F}}
\def\cJ{{\mathcal J}}

\def\cL{{\mathcal L}}
\def\cN{{\mathcal N}}

\def\cO{{\mathcal O}}
\def\cP{{\mathcal P}}

\def\cT{{\mathcal T}}

\def\qqq{\,,\,~\forall}

\newcommand{\ie}{{\it i.e.\/}\ }

\newcommand{\opcit}{{\it op.cit.\/}\ }

\def\Hom {{\mbox{Hom}}}

\def\ffp{\mathfrak{p}}

\def\sss{{\mathbb S}}
\def\gop{{\Gamma^{\rm op}}}

\def\sspec{{\pmb{\mathfrak{Spec}}}}
\def\rroot{{\pmb{\mathfrak{Root}}}}

\def\Se{\frak{ Sets}}
\def\fin{\frak{ Fin}}

\def\An{\mathfrak{ Ring}}

\def\nt{\N^{\times}}
\def\wnt{{\widehat{\N^{\times}}}}

\def\Ses{{\Se_*}}

\def\ffp{\mathfrak{p}}

\def\Mo{\mathfrak{Mo}}
\def\An{\mathfrak{ Ring}}

\def\Mr{\mathfrak{ MR}}

\def\salg{\mathfrak{ S}}

\def\Se{\mathfrak{ Sets}}

\def\smod{{\sss-{\rm Mod}}}

\begin{document}

\title{On Absolute Algebraic Geometry\\ the affine case}
\author{Alain Connes}
\email{alain@connes.org}
\address{College de France, I.H.E.S. and Ohio State University}
\author{Caterina Consani}
\email{kc@math.jhu.edu}
\address{Department of Mathematics, The Johns Hopkins
University\newline Baltimore, MD 21218 USA}
%

\begin{abstract}
 
We develop algebraic geometry for general Segal's $\Gamma$-rings and show that this new theory  unifies  two approaches we had considered earlier on (for a geometry under $\Spec \Z$). The starting observation is that the category obtained by gluing together the category  of commutative rings and that  of pointed commutative monoids, that we used in \cite{compositio} to define $\F_1$-schemes, is naturally a full subcategory of the category  of Segal's $\Gamma$-rings (equivalently of $\sss$-algebras). In this paper we develop the affine case of this general algebraic geometry: one distinctive feature is that the spectrum $\sspec(A)$ of an $\sss$-algebra is in general a Grothendieck site  rather than a point set endowed with a topology. Two striking features of this new geometry are that it is the natural domain for cyclic homology and for homological algebra as in \cite{MM}, and that new operations,  which do not make sense in ordinary  algebraic geometry, are here available. For instance, in this new context,  the quotient of a ring by a multiplicative subgroup is still an $\sss$-algebra to which our general theory applies. Thus  the adele class space gives rise naturally to an $\sss$-algebra. Finally, we show that our theory is not a special case of the T\"oen-Vaqui\' e general theory \cite{TV} of algebraic geometry under $\Spec \Z$.
\end{abstract}

\maketitle

\vspace*{6pt}\tableofcontents  
\vspace{0.1in}
\section{Introduction}
The theory of Segal's $\Gamma$-rings provides a natural extension of the theory of rings (and of semirings): in \cite{CCprel} we used this fact to extend the structure sheaf of the affine scheme $\Spec \Z$ to the Arakelov compactification $\overline{ \Spec \Z}$. While the idea implemented in this construction is the same as Durov's \cite{Durov}, the clear advantage of working with $\Gamma$-rings is  that this category forms the natural groundwork where cyclic homology is rooted. This means that de Rham theory is here naturally available. In fact, much more holds,  since the simplicial version of modules over $\Gamma$-rings (\ie the natural set-up for homological algebra in this context) forms the core of the local structure of algebraic $K$-theory \cite{DGM}. Moreover,  in \cite{CCgromov} we showed that the new $\Gamma$-ring arising as the stalk of the structure sheaf of $\overline{ \Spec \Z}$ at the archimedean place, is intimately related to hyperbolic geometry and the Gromov norm.  In \opcit we also explained that the theory of  $\sss$-modules ($\sss$ is the smallest $\Gamma$-ring corresponding to the sphere spectrum) can be thought of as the Kapranov-Smirnov theory of vector spaces over $\F_1$ \cite{KS}, provided one works in the topos of presheaves on the Segal category $\Gamma$, rather than in the topos of sets. \newline
In the present paper we clarify the relation between our work as in \cite{CCprel,CCgromov} and our previous developments on the theory of schemes over $\F_1$, originally based on some ideas of C. Soul\' e \cite{Soule}. 
More precisely, in \cite{compositio} we promoted, after the initial approach in \cite{ak}, a theory   of schemes (of finite type) over $\F_1$  using the category $\Mr$ obtained by gluing together the category $\Mo$ of commutative pointed monoids \cite{KOW, Kato, deit, TV} with the category  $\An$ of commutative rings. This  process implements a natural pair of adjoint functors relating $\Mo$ to $\An$ and follows an idea we learnt from P.~Cartier. Using the   category $\Mr$ we determined  in \cite{compositio}  a natural notion of a variety (and also of a scheme) $\mathcal X$ over $\F_1$, as a covariant functor $\mathcal X:\Mr \longrightarrow  \Se$ to the category of sets. Such a functor determines a scheme (of finite type) over $\F_1$ if it also fulfills  three properties here recalled in the Appendix  (Definition~\ref{defnfonesch}).\newline
It is well known that the category $\salg$ of $\sss$-algebras (Segal's $\Gamma$-rings) contains as full subcategories the category $\An$  (through the Eilenberg-MacLane functor $H:\An\longrightarrow \salg$) and the category $\Mo$  (by the natural functor $M\mapsto \sss M$: \cite{DGM}).
The first result of the present paper is Theorem \ref{mainthm1} stating that the category $\Mr$  appears naturally as the full subcategory of the category $\salg$, with objects either in $\An$ or in $\Mo$. Thus, the a priori artificial process of glueing $\Mr=\An\cup_{\beta,\beta^*} \Mo$ is now fully justified and suggests, in view of the results of \cite{compositio}, to develop algebraic geometry directly in $\salg$. This is what we achieve in this paper, as far as the affine case is concerned, by extending to the general case of an arbitrary $\sss$-algebra $A$,  the construction of its spectrum as a topos $\sspec(A)$, endowed with a sheaf of $\sss$-algebras. This construction determines  the prime spectrum,  when $A=HR$ for a semiring $R$, while it hands   Deitmar's spectrum \cite{deit} when $A=\sss M$ for a monoid $M$. In particular, it shows that our first development in \cite{compositio} becomes a special case of this new ``absolute" algebraic geometry.\newline
The distinctive feature of this  theory is to give meaning to operations which are meaningless in ordinary algebraic geometry,  such as taking the quotient of a ring by a subgroup of its multiplicative group, or by restricting to the unit ball in  normed rings. These operations now make sense in full generality and one can analyze their effect on the associated spectra. In fact, the construction of the spectrum $\sspec(A)$ of an $\sss$-algebra $A$ yields a Grothendieck site endowed with a presheaf of $\sss$-algebras and this datum is more refined than that of the associated topos and sheaf, as exhibited in some examples: \S \ref{sectspech1}.\newline
The paper is organized as follows: in Section \ref{sectmorphisms} we investigate morphisms of $\sss$-algebras with the goal of proving Theorem \ref{mainthm1}. The proof of this theorem  is done in \S \ref{sectthm}, after recalling in \S \ref{sectrecallmr} the older construction of the category $\Mr$. The main originality of the construction of the spectrum of an $\sss$-algebra is that it is not based, unlike the case of semirings, on a topological point space but on a topos: more precisely on a Grothendieck site whose associated topos coincides, in the special case of semirings, with the topos of the prime spectrum.\newline
The first appearance of topoi is in \S \ref{sectcatcm}, where we describe a preliminary construction taking place at the level of monoids and that will be applied later on in the paper, for an arbitrary $\sss$-algebra $A$, to the specific monoid $A(1_+)$. This preliminary construction assigns to an object $M$ of $\Mo$ a small category $C(M)$ and a corresponding topos of presheaves $\rroot(M):=\widehat {C(M)}$.  In \S \ref{sectrootz} we compute the points of $\rroot(\Z)$ and exhibit their relation with the arithmetic site  \cite{CCas}. In Section \ref{sectspecc}, we construct the spectrum $\sspec(A)$ of an arbitrary $\sss$-algebra $A$, as a Grothendieck site $(C^\infty(A),\cJ(A))$. The small category $C^\infty(A)$ is defined in \S \ref{sectcatc} by localizing the category $C(A(1_+))$, using the points of the topos $\rroot(A(1_+))$. The idea behind this localization  is to encode the expected properties of the localizations of the $\sss$-algebra $A$, with respect to multiplicative subsets of $A(1_+)$, as will be later developed in \S \ref{sectlocsalg} for the construction of the structure sheaf. The main step in the construction of  $\sspec(A)$ is the selection of the Grothendieck topology $\cJ(A)$. This part is based on the notion of partition of an object of $C^\infty(A)$, as developed in \S \ref{sectpartit}, and it depends crucially on the higher levels of the $\sss$-algebra structure, where the $n$-th higher level is used to define partitions into $n$ elements.
In \S \ref{sectmulti} we use partitions to generate multi-partitions based on rooted trees, while the Grothendieck topology is constructed in \S \ref{sectgrottop}. Finally, in \S \ref{sectfunctorial} we show that the construction of $\sspec(A)$ has the expected functoriality.\newline
 Section \ref{sectexamples} is entirely devoted to get one familiar with the construction of Section \ref{sectspecc}, by treating examples of explicit computations of $\sspec(A)$. The main result is Proposition \ref{semiring1} which identifies, for a general semiring $R$, the spectrum $\sspec(HR)$ as the Zariski site on the prime spectrum of $R$. The construction of the structure sheaf on  $\sspec(A)$ for a general $\sss$-algebra is done in Section  \ref{sectstructsheaf}. We first explain in \S \ref{sectlocsalg} that the localization of rings (or semirings) and of monoids is a special case of a general operation of localization for $\sss$-algebras $A$ with respect to a multiplicative subset of $A(1_+)$. This provides a natural presheaf of $\sss$-algebras on the category $C^\infty(A)$, and in \S \ref{sectsheafoplusplus} we construct  the structure sheaf as the associated sheaf  on $\sspec(A)$. We prove that this construction determines a sheaf of $\sss$-algebras.  \S \ref{sectsheafring} checks that our general construction agrees with the standard one in the case of rings. The functoriality of the construction is shown in \S \ref{sectfunctstrucsheaf}. The proof  is not as straightforward as in the case of rings since the construction of the structure sheaf as in \S \ref{sectsheafoplusplus} required passing to the associated sheaf.\newline Finally, in \S \ref{sectquotG} we investigate the spectrum of the quotient of a general $\sss$-algebra $A$ by a subgroup of $A(1_+)$. This applies in particular to the quotient of the ring $\A_K$ of adeles of a global field $K$ by the subgroup $K^\times$, \ie  the adele class space.  In Proposition \ref{propquot} we show  that the site $\sspec(A/G)$ is the same as $\sspec(A)$ and in Remark \ref{remsheafify} that the structure sheaf on $\sspec(A/G)$ requires to pass to the associated sheaf of the natural presheaf. In Proposition \ref{propquot1} we prove that this sheafification is not needed for $A=HR$, when the ring $R$ has no zero divisors.   \newline
   In  Section \ref{sectfinalrems} we discuss the relations between the above developments and the theory of T\"oen-Vaqui\' e \cite{TV} which is a general theory of algebraic geometry  applying to any symmetric monoidal closed category that is complete and cocomplete like the category $\smod$ of $\sss$-modules (equivalently of  $\Gamma$-sets). We prove that the theory of \cite{TV}, when implemented to the category $\smod$, does not agree with ordinary algebraic geometry when applied to the simplest example of the two point space corresponding to the spectrum of the product of two fields. The main point is Lemma \ref{twofields}  showing that the condition for a Zariski covering  in the sense of 	\cite{TV}  (Definition 2.10) is not fulfilled by the candidate cover given by the two points. \newline
For completeness, we recall in Appendix \ref{sectadjfunctor} the construction due to P. Cartier of the glueing of two categories using a pair of adjoint functors. In the same appendix we  also briefly recall the notion of a scheme over $\F_1$, as developed in \cite{compositio}.

\vspace{0.3in}

\section{Morphisms of $\sss$-algebras}\label{sectmorphisms}

In this section we establish some preliminary results on morphisms of $\sss$-modules and $\sss$-algebras which are used in Section \ref{sectglue} to identify the category $\Mr$ with a full subcategory of the category of $\sss$-algebras.

\vspace{.05in}

\subsection{Morphisms of $\sss$-modules}\vspace{.05in}
 We recall that  $\sss$-modules (equivalently $\Gamma$-sets) are by definition (covariant) functors $\gop\longrightarrow\Ses$ between pointed categories where $\gop$ is   the small, full subcategory $\Gamma^{\rm op}$  of  the category $\fin_*$ of pointed finite sets, whose objects are pointed sets $k_+:=\{0,\ldots ,k\}$, for each integer $k\geq 0$ ($0$ is the base point) and with morphism the sets $\gop(k_+,m_+)=\{f: \{0,1,\ldots,k\}\to\{0,1,\ldots,m\}\mid f(0)=0\}$. 
The morphisms  in the category $\smod$ of  $\sss$-modules are natural transformations. The category $\smod$ is a closed symmetric monoidal category (see \cite{CCgromov}).

\begin{lem}\label{scat}
$(i)$~Let $X$ be an   $\sss$-module. The map that associates to $x\in \Hom_\sss(\sss,X)$  its value $x(1_+)(1)\in X(1_+)$ on $1\in 1_+$  defines a bijection of sets $\epsilon:\Hom_\sss(\sss,X)\stackrel{\sim}{\to} X(1_+)$. \newline
$(ii)$~Let $A$ be an  abelian monoid with a zero element and $HA$ the associated $\sss$-module. There is no non-trivial morphism of $\sss$-modules from $HA$ to $\sss$.	
\end{lem}
\proof 
$(i)$~Let $j\in k_+$ and $\phi_j:1_+\to k_+$ the  map 
 $\phi_j(1)=j$. By naturality of $x\in \Hom_\sss(\sss,X)$  one has the commutative diagram
\begin{gather} \raisetag{37pt} \,\hspace{60pt}
\xymatrix@C=25pt@R=25pt{
1_+ \ar[d]_{x(1_+)} \ar[r]^{\phi_j} &
k_+\ar[d]^{x(k_+)} \\
 X(1_+) \ar[r]^{X(\phi_j)}  & X(k_+)
  } \label{natu}\hspace{100pt} 
\end{gather}
It follows that $x(k_+)(j)=X(\phi_j)\left( x(1_+)(1)\right)$. Thus, the natural transformation $x$ is uniquely determined by $\epsilon(x)=x(1_+)(1)$. Let us show that $\epsilon:\Hom_\sss(\sss,X)\to X(1_+)$ is surjective. Let $a\in X(1_+)$.
Define $x_a:\sss \to X$ as follows: to each $j\in k_+$ assign the element $x_a(k_+)(j)=X(h_{k,j})(a)\in X(k_+)$, where $h_{k,j}:1_+\to k_+$ is such that 
$h_{k,j}(1)=j$.  
One needs to check that $x_a$ is a natural transformation from the identity functor to $X$. Let $f:k_+\to n_+$ be a morphism in $\gop$. One has $f\circ h_{k,j}=h_{n,f(j)}$ and the naturality follows since 
$$
x_a(n_+)(f(j))=X(h_{n,f(j)})(a)=X(f\circ h_{k,j})(a)=X(f)\circ X(h_{k,j})(a)=X(f)( x_a(k_+)(j)).
$$
$(ii)$~The wedge sum $k_+\vee \ell_+$ of two objects of $\gop$ is defined using the inclusion maps 
$$
k_+\to (k+\ell)_+, \ \ j\mapsto j, \ \ \ell_+\to (k+\ell)_+, \ \ i\mapsto k+i.
$$
Let $\psi$ be a natural transformation from $HA$ to $\sss$. Let $\phi\in HA(k_+)=A^k$. Then $(\phi,\phi)\in A^k\times A^k=A^{2k}$ defines an element $\phi'\in HA(2k_+)=A^{2k}$ that is invariant under the permutation of the two terms in the wedge sum $k_+\vee k_+$. It follows that $\psi(\phi')=*$, since the base point is the only element of $k_+\vee k_+$ invariant under the permutation of the two terms. Let then $p:k_+\vee k_+\to k_+$ be the morphism in $\gop$ which is the identity on the first factor and the base point on the second factor. Then one has $HA(p)(\phi')=\phi$ and one concludes by naturality of $\psi$ that $\psi(\phi)=*$.\endproof\vspace{.05in}

\subsection{$\sss$-algebras and monoids}\vspace{.05in}

Let $M$ be an object of $\Mo$ \ie a multiplicative monoid  with a unit $1$ and a $0$ element. We denote by $\sss M$ the $\sss$-algebra whose underlying $\sss$-module is $M\wedge \sss$ and whose multiplicative structure is given by
\begin{equation}\label{productsm}
	\sss M(k_+)\wedge \sss M(\ell_+)\to \sss M(k_+\wedge \ell_+), \ \ (u,i)\wedge (v,j)\mapsto (uv,(i,j)).
\end{equation}
For an $\sss$-algebra $A$, the product operation restricted to $A(1_+)$ defines a multiplicative monoid with a zero (given by the base point $*$). When $A$ is commutative the monoid $A(1_+)$ is an object of $\Mo$ and this defines a functor $\cB^*:\salg \longrightarrow \Mo,~\cB^*(A)=A(1_+)$ from the category $\salg$ of $\sss$-algebras to the category $\Mo$. Next proposition shows that the functor $\cB^*$ is right adjoint to the functor $\cB:\Mo \longrightarrow \salg$, $\cB(M)= \sss M$.

\begin{prop}\label{scat1}
$(i)$~Let $A$ be an   $\sss$-algebra and $M$ be an object of $\Mo$. Then the map which to a morphism of $\sss$-algebras $\phi\in \Hom_\salg(\sss M,A)$ associates its restriction to $1_+$:
$$
\cB^*(\phi)=\phi(1_+): \sss M(1_+)=M\to A(1_+)
$$
  defines a bijection of sets $\nu:\Hom_\salg(\sss M,A)\stackrel{\sim}{\to} \Hom_\Mo(M,A(1_+))$. \newline
$(ii)$~Let $R$ be a semiring. Then there is no non-trivial morphism of $\sss$-algebras from $HR$ to $\sss M$.	
\end{prop}
\proof $(i)$~The functoriality of $\cB^*:\salg \longrightarrow \Mo$ shows that $\cB^*(\phi)=\phi(1_+)$ is a morphism in $\Mo$. By Lemma \ref{scat}, the $\sss$-module morphism underlying $\phi$ is uniquely determined by $\phi(1_+)$. This shows that the map $\nu$ is injective. Let us show that it is also surjective. Let $\psi \in \Hom_\Mo(M,A(1_+))$. By Lemma \ref{scat} there exists a unique $\sss$-module morphism $\tilde \psi\in \Hom_\sss(\sss M,A)$ which extends $\psi$.  With the notations of the proof of Lemma \ref{scat} one has, for $u\in M$ and $i\in k_+$,
$$
\tilde \psi((u,i))=A(h_{k,i})\psi(u)\in A(k_+).
$$
Let us show that $\tilde \psi$ is a morphism of $\sss$-algebras.
 The naturality of the product in $A$ determines for any pair of morphisms in $\gop$, $f:X\to X'$, $g:Y\to Y'$, a commutative diagram
\begin{gather} \raisetag{37pt} \,\hspace{60pt}
\xymatrix@C=25pt@R=25pt{
A(X)\wedge A(Y) \ar[d]_{A(f)\wedge A(g)} \ar[r]^{m_A} &
A(X\wedge Y) \ar[d]^{A(f\wedge g)} \\
A(X')\wedge A(Y') \ar[r]_{m_A}  & A(X'\wedge Y')
  } \label{diag2}\hspace{100pt}
\end{gather}
Let then $X=Y=1_+$ and $X'=k_+$, $Y'=\ell_+$, while $f=h_{k,i}$, $g=h_{\ell,j}$. One then gets from \eqref{diag2} that 
$$
m_A(A(f)\wedge A(g)(\psi(u)\wedge \psi(v))=A(f\wedge g)\psi(uv).
$$
and this equality gives, using \eqref{productsm}, the required multiplicativity.\endproof 

Next corollary uses the notation $\beta^*$ as in \eqref{units} below.

\begin{cor}\label{transit} $(i)$~Let $M$ be an object of $\Mo$, $R$ an object of $\An$ and $f\in \Hom_\Mo(M,\beta^*(R))$. Then there exists a unique morphism of $\sss$-algebras $\tilde f\in \Hom_\salg(\sss M,HR)$ such that 
$$
\nu(\tilde f)=f.
$$	
$(ii)$~Let $A$ be an $\sss$-algebra. The map of sets $\nu:\Hom_\salg(\sss[T],A)\to A(1_+)$ defines a canonical bijection  $\Hom_\salg(\sss[T],A)\simeq A(1_+)$, where $\sss[T]:=\sss M$ with $M:=\{0\}\cup \{T^n\mid n\geq 0\}$ is the free monoid with a single generator.
\end{cor}
\proof $(i)$~follows from Proposition \ref{scat1} $(i)$ and the equality $HR(1_+)=R$.\newline
$(ii)$~With $M$ as in $(ii)$ one has $\Hom_\Mo(M, A(1_+))=A(1_+)$ so the result follows  from Proposition \ref{scat1} $(i)$. \endproof

\section{The category $\Mr=\An\cup_{\beta,\beta^*} \Mo$}\label{sectglue}

In this section we prove that the category $\Mr=\An\cup_{\beta,\beta^*} \Mo$ obtained, as recalled in \S \ref{sectrecallmr}, by glueing together the categories $\Mo$ and $\An$ using an appropriate pair of adjoint functors, is in fact a full subcategory of the category of $\sss$-algebras. 

\subsection{The category $\Mr$}\label{sectrecallmr}\vspace{.05in}

In \cite{compositio} we developed the theory   of schemes (of finite type) over $\F_1$  by implementing the category obtained by gluing together the category $\Mo$ of commutative monoids \cite{KOW, Kato, deit, TV} with the category  $\An$ of commutative rings. This  process uses a natural pair of adjoint functors relating $\Mo$ to $\An$ and follows an idea we learnt from P.~Cartier. The resulting category $\Mr$  defines a  framework in which the  original definition of a variety over $\F_1$ as in \cite{Soule} is applied to a covariant functor $\mathcal X:\Mr \longrightarrow \Se$ to the category of sets. Such a functor determines a scheme (of finite type) over $\F_1$ if it also fulfills the  three properties of Definition~\ref{defnfonesch}.

The category  $\Mr=\An\cup_{\beta,\beta^*} \Mo$ is obtained by applying the  construction described in the appendix (Section \ref{sectadjfunctor}) to the following pair of adjoint covariant functors $\beta$ and $\beta^*$. The functor
\begin{equation}\label{beta}
 \beta: \Mo \longrightarrow \An\,, \quad   M\mapsto \beta(M)=\Z[M]
\end{equation}
associates to a monoid $M$ the convolution ring $\Z[M]$ (the $0$ element of $M$ is sent to $0$). The adjoint functor $\beta^*$
\begin{equation}\label{units}
 \beta^*: \An \longrightarrow \Mo\quad   R\mapsto \beta^*(R)= R
\end{equation}
 associates to a ring $R$ the ring itself viewed as a multiplicative monoid (forgetful functor). The adjunction relation states that
\begin{equation}\label{adjtrel}
\Hom_\An(\beta(M), R)\cong \Hom_\Mo(M,\beta^*(R)).
\end{equation}
We apply Proposition \ref{catplus0} to construct the category  $\Mr=\An\cup_{\beta,\beta^*} \Mo$ obtained by gluing $\Mo$ and $\An$. The collection\footnote{It is not a set: we refer for details to the discussion contained in the preliminaries of \cite{demgab}} of objects of $\Mr$ is obtained as the {\em disjoint union} of the collection of objects of $\Mo$ and $\An$. For $R\in{\rm Obj}(\An)$ and $H\in{\rm Obj}(\Mo)$, one sets $\Hom_{\Mr}(R,H)=\emptyset$. On the other hand, one defines
\begin{equation}\label{morC0bis}
    \Hom_{\Mr}(H,R):=\Hom_{\An}(\beta(H), R)\cong \Hom_{\Mo}(H,\beta^*(R)).
\end{equation}
 The morphisms between objects contained in a same category are unchanged. The composition of morphisms in $\Mr$ is defined as follows. For $\phi \in \Hom_{\Mr}(H,R)$ and $\psi \in \Hom_{\Mo}(H',H)$, one  defines $\phi\circ \psi\in \Hom_{\Mr}(H',R)$ as the composite
 \begin{equation}\label{defcomp1}
 \phi\circ \beta(\psi)\in \Hom_{\An}(\beta(H'), R)=\Hom_{\Mr}(H',R).
 \end{equation}
  Using the commutativity of the diagram  \eqref{diag}, one obtains
  \begin{equation}\label{defcomp1bis}
\Phi( \phi\circ \beta(\psi))=\Phi(\phi)\circ \psi\in \Hom_{\Mo}(H',\beta^*( R)).
 \end{equation}
  Similarly, for $\theta\in \Hom_{\An}(R, R')$ one defines $\theta \circ \phi\in \Hom_{\Mr}(H,R')$ as the composite
  \begin{equation}\label{defcomp2}
 \theta \circ \phi\in \Hom_{\An}(\beta(H), R')=\Hom_{\Mr}(H,R')
 \end{equation}
 and using  again the commutativity of \eqref{diag} one obtains that
  \begin{equation}\label{defcomp2bis}
\Phi( \theta \circ \phi)=\beta^*(\theta)\circ\Phi(\phi)
\in \Hom_{\Mo}(H,\beta^*( R')).
 \end{equation}

\subsection{The category $\Mr$ as a category of $\sss$-algebras}\label{sectthm}\vspace{.05in}

The first result which motivates the present paper is part $(iii)$ of the following theorem exhibiting $\Mr$ as a full subcategory of the category of $\sss$-algebras.

\begin{thm} \label{mainthm1}  
$(i)$~The functor $\Mo \longrightarrow \salg$, $M\mapsto \sss M$ identifies $\Mo$ with a full subcategory of $\salg$.\newline
$(ii)$~The functor $\An \longrightarrow \salg$, $R\mapsto HR$ identifies $\An$ with a full subcategory of $\salg$.\newline 
$(iii)$~The above two functors extend uniquely to a fully faithful functor $\cF:\Mr\longrightarrow \salg$ such that for any object $M$  of $\Mo$, and $R$  of $\An$ and any $f\in \Hom_\Mo(M,\beta^*(R))=\Hom_\Mr(M,R)$ one has
$$
\cF(f)=\tilde f \in \Hom_\salg(\sss M,HR).
$$
\end{thm}
\proof $(i)$~Let $M,N$ be two objects of $\Mo$. By Proposition \ref{scat1} the morphisms of $\sss$-algebras from $\sss M$ to $\sss N$ correspond canonically to elements of $\Hom_\Mo(M,N)$. \newline
$(ii)$~follows from \cite{DGM} Section 2.1. \newline 
$(iii)$~By Proposition \ref{scat1} $(ii)$ there exists no morphism of $\sss$-algebras from $HR$ to $\sss M$ independently of the choice of the objects $M$  of $\Mo$, and $R$  of $\An$. By Proposition \ref{scat1} $(i)$ and Corollary \ref{transit} $(i)$ the morphisms in $\Hom_\salg(\sss M,HR)$ (with the above notations) correspond canonically to elements of $\Hom_\Mo(M,\beta^*(R))=\Hom_\Mr(M,R)$ and one checks that the composition of morphisms in $\Mr$ corresponds to the composition of morphisms in $\salg$. \endproof 

\section{The topos $\rroot(M)$ for a monoid $M$}
 
In this section we introduce  a category $C(M)$ canonically associated to a commutative pointed monoid $M$. The topos $\rroot(M)$ dual to this small category, \ie the topos of contravariant functors $C(M)\longrightarrow \Se$, plays, in Section \ref{sectspecc}, an important role  for the construction of the Grothendieck site spectrum  $\sspec(A)$ associated to an arbitrary $\sss$-algebra $A$. The objects of $C(M)$ label the denominators in the localization process.  After defining the category $C(M)$ in \S\ref{sectcatcm}, we show in \S\ref{sectrootz} that even in the simplest case of the (multiplicative) monoid $\Z$ of the integers the computation of the points of the topos $\rroot(M)$ is quite involved.

\subsection{The category $C(M)$}\label{sectcatcm}\vspace{.05in}

Let $M$ be an object of $\Mo$. Next definition introduces the small category $C(M)$ canonically associated to $M$. 
\begin{defn} \label{cm} Let $M$ be an object  of $\Mo$.\newline 
$(i)$~We let $C(M)$ be the (small) category with one object $r(f)$ for any element $f\in M$; the morphisms are defined as follows 
$$
\Hom_{C(M)}(r(f),r(g)):= \{u\in M\mid f=ug\}.
$$
The composition of morphisms is given by the product in $M$.	\newline 
$(ii)$~$\rroot(M):=\widehat{ C(M)}$ is the topos of contravariant functors $C(M)\longrightarrow \Se$  to the category of sets.
\end{defn}
 Since for $u\in \Hom_{C(M)}(r(f),r(g))$  and $v\in \Hom_{C(M)}(r(g),r(k))$ the composition law is meaningful, one has $f=ug$, $g=vk$ so that $f=uvk$ and $uv=vu\in \Hom_{C(M)}(r(f),r(k))$.

A {\em sieve} on an object $r(f)$ of $C(M)$ is by definition a collection of morphisms with codomain $r(f)$ and stable under precomposition. A morphism of $C(M)$ with codomain $r(f)$ is of the form $u\in \Hom_{C(M)}(r(fu),r(f))$ for a unique $u\in M$. We shall use the shorthand notation $\stackrel{u}{\to}r(f)$ for $u:r(fu)\to r(f)$ since the domain of the morphism $u$ with codomain $r(f)$ is automatically $r(fu)$. The map $u\mapsto\, (\stackrel{u}{\to}r(f))$ from $M$ to morphisms with codomain $r(f)$ is bijective and compatible with the product, in particular one has a uniform description of the sieves with fixed codomain independently on this codomain.
\begin{lem}\label{sieve} The map which to an ideal $I\subset M$ of the object $M$ of $\Mo$ associates the set of morphisms of $C(M)$ with codomain $r(f)$ of the form $\stackrel{u}{\to}r(f)$ for some $u\in I$, defines a bijection between the sets of the ideals of $M$ and of sieves on $r(f)$.	
\end{lem}
When $A$ is an $\sss$-algebra, there is an interesting Grothendieck topology on $C(A(1_+))$ obtained by using Lemma \ref{sieve} and the notion of partition of unity as  in Definition \ref{partit}. However, for the purpose of constructing the spectrum of an arbitrary $\sss$-algebra we first need to   ``localize" the category $C(A(1_+))$ by inverting all morphisms of the form $f^k\in \Hom_{C(M)}(r(f^n),r(f^m))$ for $n=m+k$.  The objects of this localized category are best described as the range of a natural map from the monoid $M=A(1_+)$ to the points of the topos $\rroot(M):=\widehat{ C(M)}$. Next, we show that in the simplest example of  the monoid $M=\Z$ of the integers (with operation given by multiplication) the computation of the points of the topos $\rroot(\Z)$ involves the adeles of $\Q$ and the arithmetic site \cite{CCas}.
\vspace{.05in}

 \subsection{The points of  $\rroot(\Z)$}\label{sectrootz}\vspace{.05in}

We  consider the topos $\rroot(\Z):=\widehat{ C(\Z)}$.  The objects of $C(\Z)$ are the relative integers and the morphisms are
 $$
 \Hom_{C(\Z)}(a,b):=\{n\in \Z\mid nb=a\}.
 $$
 
 For $b\neq 0$ there exists a unique morphism $f(a,b)\in \Hom_{C(\Z)}(a,b)$ if $a\in b\Z$ while otherwise there is no morphism. 
 
 For $b=0$, $\Hom_{C(\Z)}(a,b)=\emptyset$  unless $a=0$. Finally $\Hom_{C(\Z)}(0,0)=\Z$. 
 We encode these morphisms in $C(\Z)$ as $\stackrel{n}{\to}0$. \newline
  Since the opposite category $C(\Z)^{\rm op}$ has a morphism $a\to b$ exactly when $a\vert b$, one can expect that in the process of completion--that provides the points of the topos $\widehat{ C(\Z)}$ as colimits of objects of $C(\Z)^{\rm op}$-- these limits should be naturally interpreted in terms of supernatural numbers. In fact the following proposition holds
\begin{prop} \label{ptsD}The category $\cC$ of  points of the topos $\widehat{ C(\Z)}$ contains as full subcategories  the category $\cC_1$ of subgroups of $\Q$ containing $\Z$, with morphisms  given by inclusion, and the category $\cC_2$ of ordered groups isomorphic to subgroups of $(\Q,\Q_+)$ with morphisms of ordered groups. \newline
The objects of $\cC$ form the disjoint union of the objects of the two categories $\cC_j$, $j=1,2$. There is no morphism from an object of $\cC_2$ to any object of $\cC_1$ and there is a single morphism from an object of $\cC_1$ to any object of $\cC_2$.	
\end{prop}
\proof Let $F:C(\Z)\longrightarrow\Se$ be a flat functor. We let $X_a$ be the set which is the image by $F$ of the object $a$ of $C(\Z)$. We first assume that $X_0=\emptyset$. We let $F(a,b):=F(f(a,b)):X_a\to X_b$ for $a\in b\Z$. The  filtering conditions 
\begin{enumerate}
\item $I$ is non empty.
\item For any two objects $i,j$ of $I$ there exist an object $k$ and morphisms $k\to i$, $k\to j$.
\item For any two morphisms  $\alpha, \beta:i\to j$, there exists an object $k$ and a morphism $\gamma:k\to i$
such that $\alpha\circ \gamma=\beta\circ\gamma$.
\end{enumerate}
for the category $\int_{C(\Z)}\, F$ mean that:
\begin{enumerate}
\item $X_a\neq \emptyset$ for some $a\in \Z$.
\item For any $x\in X_a$, $y\in X_b$ there exist $c\in \Z$ and $z\in X_c$  such that $a\vert c$, $b\vert c$ and $F(c,a)z=x$, $F(c,b)z=y$.
\end{enumerate}
The condition $(iii)$ is automatically satisfied since (using $X_0=\emptyset$) for objects $i,j$ of  $\int_{\cD}\, F$ there is at most one morphism from $i$ to $j$, thus $\alpha, \beta:i\to j$ implies $\alpha=\beta$. The condition $(ii)$ implies, by taking $a=b$, that there is at most one element in $X_a$ for any $a\in \Z, a\neq 0$. Let $J:=\{a\in \Z^\times  \mid X_a\neq \emptyset\}$. One then has
\begin{equation}\label{pointsD}
a\in J, \ b\vert a\implies b\in J, \  \  a,b\in J\implies \exists c\in J, \ \ a\vert c, \ b\vert c.
\end{equation}
The first implication in \eqref{pointsD} follows from the existence of the morphism $F(a,b):X_{a}\to X_b$, while the second implication derives from condition $(ii)$. Thus we get $1\in J$ and that flat functors $F:C(\Z)\longrightarrow\Se$ with $F(0)=\emptyset$  correspond to subsets $J\subset \nt$ fulfilling the conditions \eqref{pointsD}.\newline A morphism of functors $F\to F'$ exists if and only if $J\subset J'$ and when it exists is unique since each of the sets involved has one element. We now use the following correspondence between subsets $J\subset \nt$ fulfilling the conditions \eqref{pointsD} and subgroups $H$ with $\Z\subset H\subset \Q$, 
\begin{equation}\label{hj}
H^J:=\bigcup_{n\in J}\frac 1n\Z\subset \Q, \qquad  J_H:=\{n\in \Z^\times\mid \frac 1n\in H\}.
\end{equation}
The conditions \eqref{pointsD} show that $H^J$ is a filtering union of subgroups and  hence, since $1\in J$, it is a subgroup of $\Q$ that contains $\Z$. Conversely, for a group $H$ with $\Z\subset H\subset \Q$, the subset $J_H\subset \Z^\times$ fulfills  the conditions \eqref{pointsD}. Moreover, the maps $J\mapsto H^J$ and $H\mapsto J_H$ are inverse of each other. This proves that the  full subcategory of the category of points of $\widehat{ C(\Z)}$ whose associated flat functors fulfill $F(0)= \emptyset$  is the category $\cC_1$.

Let us now assume that $X_0=F(0)\neq \emptyset$. We first show that each $X_a$, for $a\neq 0$, contains exactly one element. It is non-empty since $F(0,a):X_0\to X_a$. Let $x,y\in X_a$, then using $(ii)$ and the existence of a unique morphism in $C(\Z)$ with codomain $a$ one derives $x=y$. We can then restrict $F$ to the full subcategory  of $C(\Z)$ with the single object $0$.   By applying Theorem 2.4 in \cite{CCas} one derives that  the category $\cC_2$ of points of this full subcategory of $C(\Z)$ is isomorphic to the category of totally ordered subgroups of $\Q$ with morphisms of ordered groups. There is no morphism from an object of $\cC_2$ to any object of $\cC_1$ since the corresponding natural transformation of (flat) functors cannot be defined on the object $0$ of $C(\Z)$.  Finally, we  show that there is  a single morphism from an object  of $\cC_1$ to any object  of $\cC_2$. Let $F_j$ be the corresponding (flat) functors. Then, whenever $F_1(a)\neq \emptyset$ one gets that $F_2(a)$ has a single element and is thus the final object of $\Se$, hence the existence and uniqueness of the natural transformation $F_1\to F_2$. 	 \endproof

It follows from the results of \cite{CCas} that the non-trivial subgroups of $\Q$ are parametrized by the quotient space $\A_f/\hat\Z^*$ of the finite adeles of $\Q$, while the subgroups of $\Q$ containing $\Z$ are parametrized by the subset $\hat \Z/\hat\Z^*\subset \A_f/\hat\Z^*$. This subset surjects onto the quotient $\Q_+^\times\backslash\A_f/\hat\Z^*$ which parametrizes the points of $\wnt$.

We denote by $\cD$ the full subcategory of $C(\Z)$ whose objects are the positive elements $a>0$. By construction one has 
$$
\Hom_\cD(a,b)=\begin{cases} \{f(a,b)\} &\mbox{if } b\vert a \\ 
\emptyset  & \mbox{otherwise. }\end{cases}   
$$
The proof of Proposition \ref{ptsD} then shows that the points of the topos $\widehat \cD$, dual to $\cD$, form the category $\cC_1$. 
One has a natural geometric morphism from the topos $\widehat \cD$ to the topos $\wnt$ associated to the functor
$$
\rho:\cD\longrightarrow \nt, \  \rho(a) = \bullet, \ \rho(f(a,b)):=a/b.
$$
By applying Theorem VII, 2.2 in  \cite{MM}, a functor $\phi:\cC\longrightarrow \cC'$ of small categories induces a geometric morphism $\hat \phi: \hat \cC \to \hat \cC '$  whose inverse image $\hat \phi^*$ takes a presheaf $\cF$ on $\cC'$ to the composite $\cF\circ \phi^{\rm op}$. To determine the action of the geometric morphism $\hat \phi$ on a point $\ffp$ of the topos $\hat \cC$, we consider the inverse image functor 
$$
(\hat \phi \circ \ffp)^*= \ffp^*\circ \hat \phi^*.
$$
One then identifies the flat functor $F':\cC '\longrightarrow \Se$ associated to the image $\ffp'$ of $\ffp$ using the Yoneda embedding $y:\cC '\longrightarrow {\rm Sheaves}(\hat \cC ')$, and derives 
$$
F'=(\hat \phi \circ \ffp)^*\circ y=\ffp^*\circ \hat \phi^*\circ y=\ffp^*\circ(y\circ \phi^{\rm op}).
$$
The map  $\rho$ provides  a topos theoretic meaning of the quotient map from rational adeles to adele classes.
\begin{prop} \label{ptsD1}The geometric morphism $\hat \rho: \widehat \cD\longrightarrow \wnt$  associated to the functor $\rho$ maps the category of points of the topos $\widehat \cD$ to the category of points of $\wnt$ as the restriction to $\hat \Z/\hat\Z^*\subset \A_f/\hat\Z^*$ of the quotient map
\begin{equation}\label{pointsmap1}
 \A_f/\hat\Z^* \to \Q_+^\times\backslash\A_f/\hat\Z^*
\end{equation}
\end{prop}
\proof
Given a point $\ffp$ of the topos $\widehat \cD$, let $J\subset \nt$ be the associated subset. Then the pullback $\ffp^*$ of the geometric morphism $\ffp$ associates to any contravariant functor $Z:\cD\longrightarrow\Se$ the set 
\begin{equation}\label{pointsmap}
\left(\coprod_{a\in J}Z_a\right)/\sim~ =\varinjlim_{a\in J} Z_a
\end{equation}
which is the colimit of the filtering diagram of  sets $Z_a\to Z_b$ for $a\vert b$ (remember that $Z$ is contravariant). The Yoneda embedding is given by the  contravariant functor $y(\bullet):\nt\longrightarrow \Se$ which associates to the singleton  $\bullet$ the set $\nt=\Hom_{\nt}(\bullet,\bullet)$ on which $\nt$ acts by multiplication. Once composed with $\rho$ it determines  the following contravariant  functor 
\begin{equation}\label{pointsmap2}
Z=y(\bullet)\circ \rho:\cD\longrightarrow\Se, \  \  Z(a)=\nt, \ \  Z(f(a,b))=y(\bullet)(a/b).
\end{equation}
Thus, the flat functor $F:\nt\longrightarrow\Se$ associated to the image of the point $\ffp$ by the geometric morphism $\hat\rho$ is given by the action of $\nt$ on the filtering colimit \eqref{pointsmap} applied to $Z$ of \eqref{pointsmap2}. This colimit gives 
$$
\varinjlim_{a\in J}\nt=\bigcup_{a\in J}\frac 1a\Z_+=(H_J)_+.
$$
This shows that at the level of the subgroups of $\Q$ the map of points associates to $H$, $\Z\subset H\subset \Q$ the group $H$ viewed as an abstract ordered group. This corresponds to the map \eqref{pointsmap1}. \endproof

\section{The spectrum $\sspec(A)$ of an $\sss$-algebra}\label{sectspecc}
  
This section provides the basis to develop  algebraic geometry over $\sss$. The main step, in the affine case, is to associate to an $\sss$-algebra $A$ its spectrum $\sspec(A)$ understood as the pair of a Grothendieck site and a structure presheaf of $\sss$-algebras, whose associated sheaf is a sheaf of $\sss$-algebras.  When $A=HR$ for an ordinary commutative ring $R$, $\sspec(A)$  reproduces the classical affine spectrum $\Spec R$. In this section we define the  Grothendieck site as a pair $(C^\infty(A),\cJ(A))$ of a small category $C^\infty(A)$ and a Grothendieck topology $\cJ$ on $C^\infty(A)$. The definitions of the structure presheaf and the associated sheaf are given in Section \ref{sectstructsheaf}.   The natural structure on $A(1_+)$  of a multiplicative pointed monoid (an object of $\Mo$)  provides, by the theory developed in \S \ref{sectcatcm}, a small category $C(A(1_+))$ which only depends on  $A(1_+)$ as an object of $\Mo$. The role of the category $C^\infty(A)$, obtained from $C(A(1_+))$ by applying a localization process (turning the natural morphisms between powers of the same element into isomorphisms), is to label the morphisms of localization which generalize, in the context of $\sss$-algebras, the localization morphisms $R\to R_f$  at  multiplicative subsets generated by elements $f\in R$.  This construction is described in \S\ref{sectcatc} where we also show that the category $C^\infty(A)$ has all finite limits. The definition of the Grothendieck topology $\cJ$ on $C^\infty(A)$ involves all 
the higher levels in the structure   of the $\sss$-algebra $A$ which serve as a substitute for the missing additive structure on the multiplicative monoid $A(1_+)$.  Here, the key notion  is that of a partition of an object of $C^\infty(A)$ into $n$ pieces and this process involves the $n$-th level of  the $\sss$-algebra $A$. This part is developed in \S\ref{sectpartit} and refined in \S\ref{sectmulti}.  In \S\ref{sectgrottop} we prove that one obtains in this way a Grothendieck topology $\cJ$ on $C^\infty(A)$. Lemma \ref{semirings} provides a useful tool, used later in the paper, to compare our geometric construction with the standard one in the context of semirings. Finally, in \S\ref{sectfunctorial} we  prove  Theorem \ref{functori}, stating that the Grothendieck site $(C^\infty(A),\cJ(A))$ depends functorially (and contravariantly) on the $\sss$-algebra $A$.

\subsection{The category $C^\infty(A)$}\label{sectcatc}\vspace{.05in}

Given an $\sss$-algebra $A$, the product 
$$
m_A:A(X)\wedge A(Y)\to A(X\wedge Y)
$$
together with the canonical isomorphism $1_+\wedge Y\simeq Y$ provides a meaning to the expression $uy\in A(Y)$, for $u\in A(1_+)$ and $y\in A(Y)$. We shall use freely this notation below. One has by \eqref{diag2}
\begin{equation}\label{covariance}
A(f)(ux)=uA(f)(x), \qqq f\in \Hom_\gop(X,Y), \ x\in A(X), u\in A(1_+).
\end{equation}
Next we construct, starting with the category $C(A(1_+))$, a  small category $C^\infty(A)$ which is essentially a localization of $C(A(1_+))$, turning the natural morphisms between powers of the same element into isomorphisms. The first step in the construction consists of associating to an element $b\in A(1_+)$ a point $b^\infty$ of the topos $\widehat{ C(A(1_+))}$ dual of the small category $C(A(1_+))$. The point  $b^\infty$ is  defined as the colimit of the sequence of morphisms in the opposite category $C(A(1_+))^{op} $
$$
b\to b^2 \to b^3\to \cdots 
$$ 
The understanding here is that, in general, the points of the topos which is the dual of a small category $\cC$ are obtained as filtering colimits of objects of the opposite category $\cC^{op}$. The  flat functor $F_b: C(A(1_+))\longrightarrow \Se$ associated to $b^\infty$ is defined as 
\begin{equation}\label{flatF}
F_b(a):=\varinjlim_n \Hom_{C(A(1_+))}(b^n,a).
\end{equation}
One should think of $F_b(a)$ as the set of multiplicative inverses of $a$ in the localized $\sss$-algebra of $A$ with respect to the multiplicative monoid  $M(b) :=\{b^n\mid n\in \N, n>0\}$. The proof of the first statement of the next proposition uses implicitly the uniqueness of the multiplicative inverse, when it exists.
\begin{prop}\label{atmostone} Let $A$ be an $\sss$-algebra. \newline 
$(i)$~For any $a,b\in A(1_+)$ the set $F_b(a)$ contains at most one element.\newline 
$(ii)$~One has 
	\begin{equation}\label{flatF1}
F_b(a)\neq \emptyset \iff 	\exists n\in \N, \  b^n \in aA(1_+).
\end{equation}
$(iii)$~For any  $b,c\in  A(1_+)$ the existence of a morphism of functors $F_b\to F_c$ is equivalent to 
\begin{equation}\label{flatF2}
	\exists n\in \N, \  c^n \in bA(1_+).
\end{equation}
Moreover if \eqref{flatF2} is fulfilled there exists a single morphism of functors $F_b\to F_c$.
\end{prop}
\proof $(i)$~By construction one has
$$
\Hom_{C(A(1_+))}(b^n,a)=\{c\in A(1_+)\mid ac=b^n\}
$$
and the maps involved in the (filtering) colimit are, for $n,k>0$,
$$
\beta_{n+k,n}:\Hom_{C(A(1_+))}(b^n,a)\to  \Hom_{C(A(1_+))}(b^{n+k},a), \ \ 
\beta_{n+k,n}(c):=b^k c.
$$
Let $c,c'\in \Hom_{C(A(1_+))}(b^n,a)$: we  show then that $\beta_{2n,n}(c)=\beta_{2n,n}(c')$. This equality follows from
$$
b^nc'=acc'=ac'c=b^n c.
$$
Thus the filtering colimit \eqref{flatF} contains at most one element. \newline 
$(ii)$~If $b^n \in aA(1_+)$ then $ac=b^n$ for some $c\in A(1_+)$ and the sequence $x_{n+k}=\beta_{n+k,n}(c)=b^k c$ defines an element of $F_b(a)$. Conversely, if $\Hom_{C(A(1_+))}(b^n,a)=\emptyset$ for all $n$, then $F_b(a)= \emptyset$. \newline 
$(iii)$~Let $\psi:F_b\to F_c$ be a morphism of functors. Then, since $F_b(b)\neq \emptyset$ (by $(ii)$), one has $F_c(b)\neq \emptyset$ and by $(ii)$ one has \eqref{flatF2}. Conversely, assume that \eqref{flatF2} holds, \ie $c^n=bu$ for some $n>0$ and $u\in A(1_+)$. We prove that 
$$F_b(a)\neq \emptyset \Rightarrow F_c(a)\neq \emptyset.
$$
Using $(ii)$ one gets $b^k=av$ for some $k>0$ and $v\in A(1_+)$. Thus $c^{nk}=b^k u^k=avu^k\in aA(1_+)$ and by $(ii)$ one has $ F_c(a)\neq \emptyset$. Finally by $(i)$ one gets the existence and uniqueness of  a morphism of functors $\psi:F_b\to F_c$.\endproof

\begin{defn}\label{catopen} Let $A$ be an $\sss$-algebra. We let $C^\infty(A)$ be the opposite of the full subcategory of the category of points of $\widehat{C(A(1_+))}$ of the form $b^\infty$, for $b\in A(1_+)$. 	
 \end{defn}
 
Next proposition determines $C^\infty(A)$, when $A=HR$ for a commutative ring $R$.

\begin{prop}\label{hartshorne} Let $R$ be a commutative unital ring, $\Spec R$ its prime spectrum. \newline 
$(i)$~For any two elements $b,c\in R=HR(1_+)$ one has
$$
\exists F_b\to F_c \iff V(b)\subset V(c)
$$
where $V(b):=\{\wp\in \Spec R\mid b\in \wp\}$.	
\newline 
$(ii)$~The map $D(b)\mapsto F_b$ defines an isomorphism of the category of open sets of $\Spec R$ of the form $D(b)$, $b\in R$, and morphisms given by inclusions, with the category $C^\infty(HR)$. \end{prop}
\proof $(i)$~It is well known (see \cite{Hart}, Lemma II, 2.1. (c)) that 
$$
\sqrt {bR}\supset \sqrt {cR} \iff V(b)\subset V(c).
$$
One has 
$
\sqrt {bR}\supset \sqrt {cR}\iff \exists n\in \N, \  c^n \in bR 
$.
Hence one derives by Proposition \ref{atmostone}  $(iii)$, 
$$
 V(b)\subset V(c)\iff \exists F_b\to F_c.
$$
$(ii)$~follows from $(i)$ and the uniqueness statement of Proposition \ref{atmostone}  $(iii)$, together with the equality $D(b):=V(b)^c$, which defines the open sets $D(b)\subset \Spec R$. 
\endproof 

Note in particular that if $c=0$ one always has a morphism of functors $F_b\to F_c$, for any $b\in A(1_+)$ since $F_c(a)\neq \emptyset$ for all $a\in A(1_+)$. This corresponds in algebraic geometry to the fact that the empty open set is contained in any other open set ($V(b)\subset V(0), \forall  b\Rightarrow D(0)\subset D(b), \forall b$). 
 \begin{lem}\label{finitelimits} Let $A$ be an $\sss$-algebra.
 The category $C^\infty(A)$ has all finite limits. 
 \end{lem}
 \proof  First we show that $1^\infty$ is a terminal object in $C^\infty(A)$. By Definition \ref{catopen}, this amounts to show that for any $b\in A(1_+)$ there exists a single morphism $F_1\to F_b$. This follows from  Proposition \ref{atmostone}  $(iii)$. In fact the same proposition also shows that there exists a morphism $b^\infty\to c^\infty$ in $C^\infty(A)$ if and only if there exists a morphism from some power of $b$ to some power of $c$ in the category $C(A(1_+))$. Moreover, if a morphism  $b^\infty\to c^\infty$  exists it is unique. In particular note that $(b^k)^\infty=b^\infty$ for any $k>0$. Also, note that for $f=a^\infty$, $g=b^\infty$ two objects of $C^\infty(A)$, the object $(ab)^\infty$ only depends on $f$ an $g$. In fact, as we show now, it corresponds to the pullback of the morphisms to the terminal object.  More generally, for $\phi\in \Hom_{C^\infty(A)}(a^\infty,c^\infty)$ and $\psi\in \Hom_{C^\infty(A)}(b^\infty,c^\infty)$, we prove that their pullback is given by $(ab)^\infty$ and the  unique morphisms
 $
 \phi':(ab)^\infty\to b^\infty, ~ \psi':(ab)^\infty\to a^\infty
 $. 
  By Proposition \ref{atmostone}  $(iii)$ these morphisms exist and are unique. Let then $x^\infty$ be an object of $C^\infty(A)$ and morphisms 
  $
  \alpha: x^\infty\to  a^\infty, \  \beta: x^\infty\to  b^\infty
  $
  such that $\phi\circ \alpha=\psi\circ \beta$.   By Proposition \ref{atmostone}  $(iii)$ there 
   exists a power $x^n$ divisible by $a$, and  $x^m$  divisible by $b$. Thus $x^{n+m}$ is divisible by $ab$ and one obtains a (unique) morphism  $x^\infty\to (ab)^\infty$. The uniqueness of morphisms when they exist thus shows that $(ab)^\infty$ is the pullback.\endproof  \vspace{.05in}
   
   \subsection{Partitions in $C^\infty(A)$}\label{sectpartit}\vspace{.05in}
   
   In classical algebra the operation of ``sum" of two elements
  can be iterated  by applying the associative rule to obtain the sum of $n$ elements. Then one defines
a partition of unity as $n$ elements whose sum is equal to 1. In the context
of $\sss$-algebras the sum of $n$ elements cannot be well defined using only
the first two levels of an $\sss$-algebra $A$: at the best the first two levels  give  a meaning to the statement $1\in a+b$ for $a,b\in A(1_+)$. In general we know that $a+b$, when it is defined,  is a hypersum. Thus one  uses the level $n$ (\ie $A(n_+)$) 
 to give a meaning to the statement $1\in\sum_{j=1}^n a_j$   and this defines 
  a partition of unity in $n$ pieces. Lemma \ref{Gpart2}   asserts the crucial statement that one can take the product of two
partitions of unity. If the first is a partition in $n$ pieces and the second in $m$
pieces, then the product is a partition in $nm$ pieces. This proves that 
if one tried to restrict to  level two this construction would not be 
stable by product.  

   Let  $A$ be an $\sss$-algebra. We use the $\sss$-module structure of $A$ to define partitions of objects of $C^\infty(A)$. First we introduce some notations. Given a finite set $F$, we define specific elements $\delta_j$, $j\in F$, and $\Sigma_F$ of $\Hom_\gop(F_+,1_+)$ as follows
$$
\delta_j(i):=1\iff j=i, \ \ \Sigma_F(i):=1 \qqq i\in F.
$$
   \begin{defn}\label{partit} Let $A$ be an $\sss$-algebra, $f$ an object of $C^\infty(A)$. A {\em partition} of $f$ is a finite collection of morphisms $f_j\to f$, $j=1, \ldots ,n$, in $C^\infty(A)$ such that there exists $\xi\in A(n_+)$ with   $(A(\delta_j)(\xi))^\infty=f_j$, $\forall j$, and $(A(\Sigma) \xi)^\infty=f$.   	\end{defn}
   
  One has the following multiplicativity property
  
\begin{lem}\label{Gpart2} Let $A$ be an $\sss$-algebra. Let $(f_i)_{i\in F}$ and $(g_j)_{j\in G}$ be partitions of objects $u,v$ of $C^\infty(A)$. Then the family $(f_ig_j)_{(i,j)\in F\times G}$ is a partition of $uv$.
\end{lem}
\proof Let $\xi \in A(F_+)$ (resp. $\eta \in A(G_+)$)  such that $A(\delta_i)(\xi)=a_i$, $a_i^\infty=f_i$, $\forall i\in F$, and $A(\Sigma) \xi=s$, $s^\infty=u$ (resp. $A(\delta_j)(\eta)=b_j$, $b_j^\infty=g_j$, $\forall j\in G$, and $A(\Sigma) \eta=t$, $t^\infty=v$). Consider the diagram \eqref{diag2} in the case $X=F_+$, $Y=G_+$, $X'=1_+$, $Y'=1_+$:
\begin{gather} \raisetag{37pt} \,\hspace{60pt}
\xymatrix@C=25pt@R=25pt{
A(F_+)\wedge A(G_+) \ar[d]_{A(\alpha)\wedge A(\beta)} \ar[r]^{m_A} &
A(F_+\wedge G_+) \ar[d]^{A(\alpha\wedge \beta)} \\
A(1_+)\wedge A(1_+) \ar[r]_{m_A}  & A(1_+\wedge 1_+=1_+)
  } \label{diag2bis}\hspace{100pt}
\end{gather}
Taking $\alpha=\delta_i$, $\beta=\delta_j$ one obtains $\alpha\wedge \beta=\delta_{(i,j)}$ and this gives the equality, with $\zeta:=m_A(\xi \wedge \eta)$
$$
A(\delta_{(i,j)})(\zeta)= a_ib_j.
$$
Taking $\alpha=\Sigma_F$, $\beta=\Sigma_G$ one obtains $\alpha\wedge \beta=\Sigma_{F\times G}$ and this gives the equality $A(\Sigma_{F\times G}) \zeta=st$. Finally one has $( a_ib_j)^\infty=f_ig_j$ and $(st)^\infty=uv$.	\endproof
\vspace{.05in}

\subsection{Rooted trees and multi-partitions}\label{sectmulti}\vspace{.05in}

 Next we construct, for an arbitrary $\sss$-algebra $A$, a Grothendieck topology $\cJ(A)$ on $C^\infty(A)$.   For  a commutative unital ring $R$, we shall show  that the corresponding Grothendieck topology on $A=HR$ coincides with the topology of $\Spec R$ in the standard sense.\newline
 Iterating the process of partitioning  an object of $C^\infty(A)$ is best labelled by a combinatorial datum called a rooted tree.
 We start by recalling  the terminology pertaining to the notion of rooted tree. A rooted tree $(T,r)$ is a connected graph $T$ which is a tree and is  pointed  by the choice of a vertex $r$ called the root. 
 \begin{figure}	\begin{center}
\includegraphics[scale=0.7]{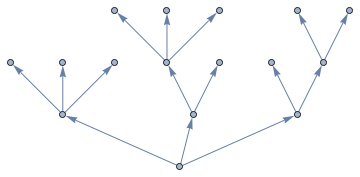}
\end{center}
\caption{Rooted tree. \label{rooted} }
\end{figure}
We shall consider a rooted tree $(T,r)$ as a small category with objects the vertices of the tree and morphisms generated by the identity morphisms and the oriented edges of the graph, for the orientation of edges as shown in Figure \ref{rooted}. The {\em height} of a  rooted tree is the maximal length of paths.  Each internal vertex $v$ has successors  which form a subset $S(v)$ of the set of vertices. The vertices with no successor are called external vertices and their collection is denoted $\cE(T)$. 
\begin{defn}\label{rootedpart} Let $A$ be an $\sss$-algebra, and $(T,r)$ a rooted tree. A contravariant functor $F:(T,r)\longrightarrow C^\infty(A)$ is called a {\em multi-partition} of $F(r)$ if for each internal vertex $v$ the collection $(F(u))_{u\in S(v)}$ forms a partition of $F(v)$ in the sense of Definition \ref{partit}.	
\end{defn}\vspace{.05in}

 \subsection{The Grothendieck topology on $C^\infty(A)$}\label{sectgrottop}\vspace{.05in}
 
  Given a small category $\fC$, an object $c$ of $\fC$ and a set $X$ of morphisms with codomain $c$, we denote by ${\rm Sieve}(X)$ the sieve generated by $X$.  We consider the notion of covering in the category $C^\infty(A)$ given by the function $K$ which associates to each object $c$ the collection $K(c)$ of families $X$ of morphisms with codomain $c $ for which there exists a multi-partition $ F:(T,r)\longrightarrow C^\infty(A)$ with $F(r)=c$, such that for all external vertices $e$ of $T$
 \begin{equation}\label{covfamil}
 	F(r\to e)\in {\rm Sieve}(X).
 \end{equation}
Note that one has 
$$
{\rm Sieve}(X)=\cup_{x\in X} {\rm Sieve}(x),  \qquad {\rm Sieve}(x)=\{y^\infty x\mid y\in A(1_+)\}.
$$
 
 \begin{lem}\label{finitelimits} Let $A$ be an $\sss$-algebra.
The function $c\mapsto K(c)$ defines a basis for a Grothendieck topology $\cJ(A)$ on $C^\infty(A)$.
 \end{lem}
 \proof     We first check that the function that associates to an object $c$ of $C^\infty(A)$ the collection $K(c)$ of families of morphisms of \eqref{covfamil} fulfills the three conditions of Definition III. 2 of \cite{MM}. The first condition states that the family with a single element $c\to c$ belongs to $K(c)$. This follows taking the trivial partition of $c$. The second condition is the stability by pullback, which follows if we show  
   that for $X\in K(c)$ and $d\to c$ a morphism, the family $dX:=(xd)_{x\in X}$ belongs to $K(d)$. Indeed, one has   $cd=d$.  Let $ F:(T,r)\longrightarrow C^\infty(A)$ be a multi-partition with $F(r)=c$ such that \eqref{covfamil} holds.  Lemma \ref{Gpart2} shows that the functor $dF:(T,r)\longrightarrow C^\infty(A)$ defined as the pointwise product of $F$ with $d$ fulfills \eqref{covfamil}, with respect to $dX$ which thus belongs to $K(d)$.    
   The third condition of Definition III. 2 of \opcit states that if $X=\{f_i:c_i\to c\mid i\in I\}$ belongs to $K(c)$ and for each $i\in I$ one has a family $Y_i=\{g_{ij}:d_{ij}\to c_i\mid  j\in I_i\}$, $Y_i\in K(c_i)$, then the family $Y:=\{f_i\circ g_{ij}:d_{ij}\to c\mid i\in I, j\in I_i\}$ is in $K(c)$. By hypothesis there exists  a multi-partition $ F:(T,r)\longrightarrow C^\infty(A)$ with $F(r)=c$ such that \eqref{covfamil} holds, and for each $i\in I$ a multi-partition $ F_i:(T_i,r_i)\longrightarrow C^\infty(A)$ with $F_i(r_i)=c_i$ such that for all external vertices $\epsilon$ of $T_i$
 \begin{equation}\label{covfamili}
 	F_i(r_i\to \epsilon)\in {\rm Sieve}(Y_i).
 \end{equation}
Let $e$ be an external vertex of $T$. One 
 has, by \eqref{covfamil}, $ 	F(r\to e)\in {\rm Sieve}(X)$, so let $\psi(e)\in I$ be such that $ 	F(r\to e)\in {\rm Sieve}(c_{\psi(e)})$. We thus obtain a map $\psi:\cE(T)\to I$ from the set $\cE(T)$ of external vertices of $T$ to $I$ and a map $\alpha:\cE(T)\to A(1_+)$ such that 
 $$
 F(r\to e)=\alpha(e)^\infty c_{\psi(e)}.
 $$
  We let then $(T',r)$ be the rooted tree obtained by grafting the rooted tree $(T_{\psi(e)},r_{\psi(e)})$ at each external vertex $e$ of $T$. Then we define a multi-partition $ F':(T',r)\longrightarrow C^\infty(A)$ as follows. One has $F_i(r_i)=c_i$ for all $i\in I$ and hence for $i=\psi(e)$ one gets 
  $$
  c_{\psi(e)}=F_{\psi(e)}(r_{\psi(e)}), \  \ F(r\to e)=\alpha(e)^\infty F_{\psi(e)}(r_{\psi(e)})
  $$
which shows that the functors $\alpha(e)^\infty F_{\psi(e)}$ match with $F$ and define a single functor $ F':(T',r)\longrightarrow C^\infty(A)$. It is a multi-partition since the condition of Definition \ref{rootedpart} is checked locally for each internal vertex. Each external vertex $e'$ of $T'$ is an external vertex of a $(T_{\psi(e)},r_{\psi(e)})$. By \eqref{covfamili} one has $F'(r\to e')\in {\rm Sieve}(Y)$ thus one gets $Y\in K(c)$ as required.\endproof 

\begin{defn}\label{defnspecc} Let $A$ be an $\sss$-algebra. We define the spectrum $\sspec(A)$ as the Grothendieck site $(C^\infty(A),\cJ(A))$.	
\end{defn}

Note that the information contained in the site is more precise than that of the associated topos: 
in \S\ref{sectspech1} we shall give an example where this nuance plays a role. 

Next lemma describes the Grothendieck topology $\cJ(A)$ for $A=HR$, where $R$ is a semiring.

\begin{lem}\label{semirings} Let $R$ be a semiring and $A=HR$ the associated $\sss$-algebra.\newline
$(i)$~Let $f$ be an object of $C^\infty(A)$,  $f_j\to f$, $j=1, \ldots ,n$ a partition of $f$ in the sense of Definition \ref{partit}. Let $x_j, x\in R$ such that $x_j^\infty=f_j$ for $j=1, \ldots ,n$ and $x^\infty=f$. Then there exist elements $a_j\in R$ such that $\sum_{j=1}^n a_jx_j=x^k$, for some power $x^k$ of $x$. \newline
$(ii)$~Let $F:(T,r)\longrightarrow C^\infty(A)$ be a  multi-partition of $F(r)$, $x,x_e\in R$  such that $x^\infty =F(r)$ and $x_e^\infty=F(e)$, for all $e\in \cE(T)$. Then there exist elements $a_e\in R$ such that $\sum a_ex_e=x^k$ for some power $x^k$ of $x$.\newline
$(iii)$~A family $X$ of morphisms $x^\infty\to c^\infty$ with codomain $c^\infty$ is a covering for the  Grothendieck topology $\cJ(A)$ if and only if there exists a map $a:X \to R$ with finite support such that $\sum a(x) x$ is equal to a power of $c$.
 \end{lem}
 \proof $(i)$~Let $\xi\in HR(n_+)$ with   $(HR(\delta_j)(\xi))^\infty=f_j$, $\forall j$, and $(HR(\Sigma) \xi)^\infty=f$. One has $HR(n_+)=R^n$ and thus one gets elements  $\xi_j\in R$, $j=1,\ldots,n$, such that  $\xi_j^\infty=f_j=x_j^\infty$ for all $j$ and $(\sum_1^n \xi_j)^\infty=x^\infty$. Thus, each $x_j$ divides a power of $\xi_j$ and we let $m$ and $b_j\in R$ such that $\xi_j^m=b_jx_j$ for all $j$. Let $\xi=\sum_1^n \xi_j$. Then, there exist elements $c_j\in R$ such that 
 $$
 \xi^{mn}=(\sum_1^n \xi_j)^{mn}=\sum_1^n c_j\xi_j^m.
 $$
 One has $\sum_1^n c_j\xi_j^m=\sum_1^n c_jb_jx_j$ and since $\xi^{\infty}=x^\infty$,  $\xi^{mn}$ divides a power of $x$, \ie $x^k=a \xi^{mn}$ for some $a\in R$. Thus one gets $x^k=\sum_1^n a c_jb_jx_j=\sum a_jx_j$ as required.
 \newline
$(ii)$~Let us choose for each vertex $w$ of $T$ an element $x(w)\in R$ such that $F(w)=x(w)^\infty$ and that $x(e)=x_e$ for external vertices, while $x(v)=x$. Consider a vertex $w$ such that all its successors (\ie vertices in $S(w)$) are external. Then by hypothesis the $(F(u))_{u\in S(w)}$ forms a partition of $F(w)$ and hence by $(i)$ one can find a power $x(w)^k$ of the form $\sum_{S(w)} a(u) x(u)$. One  has $F(w)=(x(w)^k)^\infty$. Thus at the expense of raising the $x(w)$ to some power one can assume that for each vertex $t$ of $T$ one has elements $a(t,u) \in R$, $u\in S(t)$, such that
$$
x(t)=\sum_{u\in S(t)} a(t,u) x(u).
$$
 It follows then that  there exists elements $a_e\in R$ such that $\sum a_ex_e=x^k$ for some power of $x$. 
  \newline
$(iii)$~Let $f=c^\infty$, $X\in K(f)$ and $ F:(T,r)\longrightarrow C^\infty(A)$ with $F(r)=f$ such that \eqref{covfamil} holds. This shows that for each $e\in \cE(T)$ there exists an element $x(e)\in X$ and a $b(e)\in R$ such that $F(e)=(b(e)x(e))^\infty$. Thus by $(ii)$ one can find  elements $a(e)\in R$ for $e\in \cE(T)$  and $k>0$ with$$
c^k=\sum_{e\in \cE(T)} a(e)b(e)x(e)
$$
so that $X$ satisfies the condition of the Lemma. Conversely assume that there exists a map $a:X \to R$ with finite support such that $\sum a(x) x$ is equal to a power $c^k$ of $c$. Consider the partition of $c^\infty$ given by the $(a(x) xc)^\infty$. One can view it as a multi-partition of height $1$. For each external vertex one has $(a(x) xc)^\infty\in {\rm Sieve}(X)$ since $x^\infty\in X$. It follows that $X\in K(f)$. \endproof \vspace{.05in}

\subsection{Functoriality of $\sspec(A)$}\label{sectfunctorial}\vspace{.05in}

Let $A,B$ be $\sss$-algebras and $\phi:A\to B$ a morphism of $\sss$-algebras. We show that $\phi$ determines a morphism of sites (in the sense of Grothendieck) 
$$
\tilde \phi: (C^\infty(A),\cJ(A))\longrightarrow  (C^\infty(B),\cJ(B)).
$$
One first defines a covariant functor which assigns to an object $c=f^\infty$ of $C^\infty(A)$ the object $\tilde \phi(c):=\phi(f)^\infty$ of $C^\infty(B)$. This definition makes sense since 
$$
f^n \in gA(1_+)~ \Longrightarrow~ \phi(f)^n\in \phi(g)B(1_+).
$$  
Moreover, this implication shows that $\tilde \phi$ defines a covariant functor $C^\infty(A)\longrightarrow  C^\infty(B)$. This functor preserves finite products since $\phi(ab)=\phi(a)\phi(b)$, $\forall a,b\in A(1_+)$. It also preserves the terminal object $1$ as well as the equalizers (since they are all trivial). Thus  $\tilde \phi$ preserves finite limits. Let us show that it preserves covers. This means that if $S$ is a covering sieve of the object $c$ of $C^\infty(A)$, then the sieve generated by the morphisms $\tilde \phi(u)$, for $u\in S$, is a covering sieve of the object $\tilde \phi(c)$
of $C^\infty(B)$. By definition of the Grothendieck topology $\cJ(A)$, the sieve $S$ contains a family $X$ of morphisms with codomain $c$ for which there exists a multi-partition $ F:(T,r)\longrightarrow C^\infty(A)$ with $F(r)=c$ such that for all external vertices $e$ of $T$ one has \eqref{covfamil} \ie $
 	F(r\to e)\in {\rm Sieve}(X)$. Thus $S$ in fact contains all the morphisms $F(r\to e)$. The composition $G:=\tilde \phi\circ F$ is a multi-partition $ G:(T,r)\longrightarrow C^\infty(B)$ and thus the sieve generated by the morphisms $\tilde \phi(u)$ for $u\in S$ contains all the morphisms $G(r\to e)$ and is hence a covering sieve of the object $\tilde \phi(c)$
of $C^\infty(B)$. We  now apply Theorem VII, 10, 2 of \cite{MM} and obtain
\begin{thm} \label{functori}
Let $A,B$ be $\sss$-algebras and $\phi:A\to B$ a morphism of $\sss$-algebras. Then composition with $\phi$ defines a morphism of sites $\tilde \phi: (C^\infty(A),\cJ(A))\longrightarrow  (C^\infty(B),\cJ(B))$ and an associated geometric morphism of toposes $f:\sspec(B)\longrightarrow\sspec(A)$.	The direct image functor $f_*$ sends a sheaf $\cF$ on $\sspec(B)$ to the composite $f_*(\cF)=\cF\circ \tilde \phi$.\end{thm}

\section{Examples of $\sspec(A)$} \label{sectexamples}

\subsection{The case of rings}\vspace{.05in}

Next proposition can be deduced directly from Lemma \ref{semirings} but we give an independent and more geometric proof.
\begin{prop}\label{ordinary} Let $R$ be a commutative unital ring. The site   $(C^\infty(HR),\cJ(HR))$ is the same as the site associated to the topological space $\Spec(R)$ and the category of open sets of the form $D(f)$, $f\in R$.	
\end{prop}
\proof Let  $A=HR$. Then by Proposition \ref{hartshorne}   the category $\cC:=C^\infty(HR)$ is the same as the category of open sets of the form $D(f)\subseteq \Spec(R)$.	 Let  $f_j\to f$, $j=1, \ldots ,n$ be a partition of  an object $f\in C^\infty(A)$ as in Definition \ref{partit}. Let $\xi\in A(n_+)$, with   $(A(\delta_j)(\xi))^\infty=f_j$, $\forall j$, and $(A(\Sigma) \xi)^\infty=f $. One has $A(n_+)=R^n$ and 
$\xi=(\xi_j)_{j=1, \ldots ,n}$, moreover $(A(\delta_j)(\xi))=\xi_j$ and $A(\Sigma) \xi=\sum \xi_j$. The equality $(A(\delta_j)(\xi))^\infty=f_j$ means $D(f_j)=D(\xi_j)$, while $(A(\Sigma) \xi)^\infty=f $ means  $D(f)=D( \sum \xi_j)$. Thus if such a partition exists one has 
$
D(f)=\cup D(f_j)
$.
Hence the existence of a multi-partition  $F:(T,r)\longrightarrow C^\infty(A)$ of $f=F(r)$ implies 
$$
D(f)=\cup_{e\in \cE(T)} D(F(e)).
$$ 
This shows that a family $X\in K(c)$ is a covering of $D(c)$ in the usual sense since each open set  $D(F(e))$ of the ordinary covering associated to the multi-partition  is contained in one of the form $D(x)$ for some $x\in X$. Conversely,  if a family $X$ of objects of $C^\infty(HR)$ is a covering in the usual sense, \ie  if $\cup_{x\in X} D(x)=D(c)$, then the ideal generated by $X$ contains a power $c^n$ of $c$ and there exists a finite decomposition  $c^n=\sum x\, a(x)$, with $a(x)\in R$ which provides a partition $c_j$, ${j=1, \ldots ,n}$ of $c$ in the sense of Definition \ref{partit}. Moreover, each $c_j$ belongs to the sieve of $X$ since it is of the form $x\, a(x)$ for some $x\in X$. Thus the usual covering $X$ belongs to $K(c)$ (and the involved rooted tree is of height $1$). \endproof \vspace{.05in}

\subsection{The case of monoids} \vspace{.05in}

Let $M$ be a multiplicative, commutative, pointed  monoid (\ie an object of $\Mo$). Let $A=\sss M$ be the associated $\sss$-algebra. The spectrum $\Spec \F_M$ is the topological space whose points are the prime ideals  $I\subsetneq M$ of $M$ whose complement is stable under multiplication \cite{deit}. This spectrum always contains the prime ideal $\ffp_M=(M^\times)^c$.  The topology of $\Spec \F_M$ is defined by its closed sets which are either empty or associated to an ideal $J\subset M$ by 
$$
V(J):=\{ \ffp \in \Spec \F_M\mid J \subset \ffp\}.
$$
For any ideal $J\subset M$ one lets $\sqrt {J}:=\{a\in M\mid \exists n, \ a^n\in J \}$.
In the next part we use the notations of Proposition \ref{atmostone}.

     \begin{prop}\label{hartshorne1} Let $M$ be a multiplicative, commutative and pointed  monoid. \newline 
$(i)$~For any two elements $b,c\in M=\sss M(1_+)$ one has
$$
\exists F_b\to F_c \iff V(b)\subset V(c)
$$
where $V(b):=\{\wp\in \Spec \F_M \mid b\in \wp\}$.	
\newline 
$(ii)$~The map $D(b)\mapsto F_b$ defines an isomorphism of the category of open sets of $\Spec \F_M$ of the form $D(b)$, $b\in M$, with the category $C^\infty(\sss M)$.\newline 
$(iii)$~The site $(C^\infty(\sss M),\cJ(\sss M))$ is the same as the site associated to the space $\Spec(\F_M)$ with topology defined by the category of open sets of the form $D(f)$.	\end{prop}
 \proof $(i)$~By applying Theorem 1.1 in \cite{Gilmer},  one has 
$$
\sqrt {bM}\supseteq \sqrt {cM} \iff V(b)\subseteq V(c).
$$
Also, one has 
$$
\sqrt {bM}\supseteq \sqrt {cM}\iff \exists n\in \N, \  c^n \in bM.
$$
Hence, by Proposition \ref{atmostone}  $(iii)$, one gets 
$$
 V(b)\subseteq V(c)\iff \exists F_b\to F_c.
$$
$(ii)$~follows from $(i)$ and the uniqueness statement of Proposition \ref{atmostone}  $(iii)$, together with the equality $D(b):=V(b)^c$ which defines the open sets $D(b)$. \newline 
$(iii)$~Let $b\in M$ and $D(b):=\{\wp\in \Spec \F_M \mid b\notin \wp\}$ the associated open set. Assume $D(b)\neq \emptyset$. Let $\wp_b:=\cup_{\wp \in D(b)} \wp$: it is the largest prime ideal not containing $b$. Let then 
$$
(D(f_j))_{j\in I}, \  \  \cup D(f_j)=D(b)
$$
be an open covering of $D(b)$ of the topological space $\Spec(\F_M)$. There exists $j\in I$ such that $\wp_b\in D(f_j)$. One then has $f_j\notin \wp_b$ and it follows from the definition of $\wp_b$ that $f_j\notin \wp$, $\forall \wp \in D(b)$. This shows that $D(f_j)=D(b)$,  thus the Grothendieck topology associated to the basis of open sets of the form $D(f)$ is the chaotic topology in which every presheaf is a sheaf. Let us now show that the same holds for the Grothendieck topology $\cJ(\sss M)$. Let $f\in M$. Consider a partition  of $f$ in the sense of Definition \ref{partit} \ie a finite collection of morphisms $f_j\to f$, $j=1, \ldots ,n$, in $C^\infty(A)$ such that there exists $\xi\in A(n_+)$, with   $(A(\delta_j)(\xi))^\infty=f_j$, $\forall j$, and $(A(\Sigma) \xi)^\infty=f$. By construction of $A(n_+)=M\wedge n_+$, one has $\xi=m\wedge k$, for some $m\in M$ and $k\in n_+$. Thus $A(\delta_j)(\xi)=m$ for $j=k$ and $A(\delta_j)(\xi)=*$ otherwise. Moreover $A(\Sigma) \xi=m$. Thus such a partition of $f$ consists of a single element equal to $f$. This implies that the same holds for a multi-partition. Hence the Grothendieck topology $\cJ(\sss M)$ is the chaotic topology and agrees with that of $\Spec(\F_M)$.
\endproof 
In particular the proof of $(iii)$ in the above proposition shows that the topos $\sspec(\sss M)$ is of presheaf type. \vspace{.05in}

\subsection{The tropical case} \label{secttropical}\vspace{.05in}

Let $I\subseteq \R$ be an open bounded interval and $R=C_{\rm Conv}(I)$ be the semiring of continuous, piecewise affine with integral slope, convex functions $f:I\to \R$. In the following we consider only functions with finitely many slopes and moreover  we add in $R$ the constant function with value $-\infty$, which plays the role of the zero element of $R$. The two operations in $R$ are the pointwise sup for addition, and the pointwise sum for multiplication. We let $A:=HR$. Our first task is to determine $C^\infty(A)$. Given $f\in R$ we let $f'':=(\frac{d}{dx})^2f$ be its second derivative viewed as a distribution. It is by construction a finite sum of Dirac masses with positive integral coefficients
\begin{equation}\label{dirac}
	f''(x)=\sum n_j \, \delta_{x_j}(x)=\sum_{y\in Z(f)}n_f(y) \delta_y(x).
\end{equation}
This equality defines uniquely the finite subset (tropical zeros) $Z(f)\subset I$ and the positive integer valued function $n_f:I\to \N$ that vanishes outside $Z(f)$.

\begin{lem}\label{tropdiv} Let $f,g\in R$ be non-zero\footnote{\ie elements not equal to the constant function $-\infty$} elements. \newline
$(i)$~The function $g$ divides $f$ in $R$ if and only if $n_g\leq n_f$. \newline
$(ii)$~One has 
\begin{equation}\label{flatF3}
	\exists k\in \N, \  f^k \in gA(1_+)  \iff Z(g)\subset Z(f).
\end{equation}	
\end{lem}
\proof $(i)$~Since the product in $R$ is given by the pointwise sum, the map $f\mapsto n_f$ transforms product into sums. It follows that if $g$ divides $f$ one has $n_g\leq n_f$. Conversely if $n_g\leq n_f$, let $h$ be a function which is the double primitive of the distribution $\sum (n_f(y)-n_g(y)) \delta_y$. The choice of the double primitive is unique up to an affine function, and since such elements of $R$ are invertible, this choice only affects the equality $f=gh$  up to an invertible element. Thus one checks that $g$ divides $f$ in $R$.  \newline
$(ii)$~One has $n_{f^k}=k\, n_f$, $Z(f^k)=Z(f)$. If  $f^k \in gA(1_+)$ one has $n_g\leq k \, n_f$ from $(i)$ so that $Z(g)\subset Z(f)$. Conversely, if $Z(g)\subset Z(f)$, there exists a finite $k\in \N$ such that $n_g\leq k \, n_f$ and $(i)$ shows that $f^k \in gA(1_+)$. \endproof

\begin{prop}\label{hartshorne2} Let $R=C_{\rm Conv}(I)$. \newline 
$(i)$~For any two non-zero elements $b,c\in R=HR(1_+)$ one has
$$
\exists F_b\to F_c \iff Z(b)\subset Z(c).
$$
$(ii)$~The map $Z(b)^c\mapsto F_b$ defines an isomorphism of the category $\cF$ of open sets of $I$ which are either empty or  whose complement is finite, with the category $C^\infty(HR)$.\newline 
$(iii)$~The coverings of an object $\Omega$ for the site $(C^\infty(HR),\cJ(HR))$ are the families of objects $(\Omega_j)_{j\in K}$ of $\cF$ which are coverings of $\Omega$ in the usual sense, and such that each connected component  $C\subset \Omega$ is covered by a single element $\Omega_{j(C)}$ of the family.	\newline
$(iv)$~The category of points of the topos $\sspec(HR)$ is canonically isomorphic to the category ${\rm Conv}(I)$ of convex subsets of the interval $I$, with morphisms given by reverse inclusion.
 \end{prop} 
 \proof $(i)$~This follows from \eqref{flatF2} in Proposition \ref{atmostone} combined with \eqref{flatF3} in Lemma \ref{tropdiv}.\newline 
$(ii)$~For the zero element of $R$, \ie the function which is constant equal to $-\infty$ one defines $Z(0)=I$. By Definition \ref{catopen} the category $C^\infty(HR)$ is the opposite of the category formed by the flat functors $F_b$. Since every finite subset of $I$ appears as a $Z(f)$ using \eqref{dirac}, one gets the required isomorphism using $(i)$.\newline 
$(iii)$~Let $f$ be a non-zero element of $R=HR(1_+)$. Let $C$ be a connected component of the open set $Z(f)^c$. Let $(f_j)_{j\in K}$ be a finite family of elements of $R$ such that the pointwise supremum $\vee f_j=f$. The function $f$ restricted to $C$ is affine and there exists an index $j$ such that $f_j$ agrees with $f$ on a non-empty interval inside $C$. Then since $f_j$ is convex and $f_j\leq f$ one has $f_j(x)=f(x)$, $\forall x\in C$. Indeed, one can assume that $f(x)=0$, $\forall x\in C$ and use the convexity of $f_j$ to get that its upper-graph is bounded below by the segment $\{(x,0),x\in C\}$. This shows that a partition, in the sense of Definition \ref{partit}, of the object $f^\infty$ of $C^\infty(HR)$ contains an element $f_j$ such that $Z(f_j)\cap C=\emptyset$. By following the relevant branch of the tree, one sees that for any multi-partition of $f^\infty$ $, F:(T,r)\longrightarrow C^\infty(HR)$  there exists an external vertex $e$ of $T$ such that $Z(F(e))\cap C=\emptyset$. Hence any covering of the open set $Z(f)^c$ for the Grothendieck topology $\cJ(HR)$ fulfills the condition of $(iii)$. Conversely, let us show that if an ordinary covering $(\Omega_j)_{j\in K}$ of an object $\Omega$ fulfills the condition $(iii)$ one can find an $f\in R$ with $f^\infty=\Omega$ and elements $f_j\in R$, with $f_j^\infty=\Omega_j$, such that   $f=\vee f_j$ and hence that $(\Omega_j)_{j\in K}$ is a cover   of $\Omega$ for the Grothendieck topology $\cJ(HR)$. For each connected component $C$ of $\Omega$, let $\Omega_{j(C)}$ cover $C$ and $h_{j(C)}$ with $h_{j(C)}^\infty=\Omega_{j(C)}$ be such that it is constant equal to $0$ on $C$. Let $h\in R$ with $h^\infty=\Omega$ and for each connected component $C$ of $\Omega$ let $L_C$ be the affine function which agrees with $h$ on $C$. Then, by convexity, one has that $(h-L_c)(x)$ is strictly positive outside $C$ and bounded below by the distance $d(x,C)$. This shows that  for $n$ large enough one gets
$$
nh=\vee_C (h_{j(C)}+nL_C),
$$
which gives the required equality for $f:=nh$, $f_{j(C)}=h_{j(C)}+nL_C$, while the other $f_j$'s do not matter. \newline 
$(iv)$~A point of the topos $\sspec(HR)$  associated to the site $(C^\infty(HR),\cJ(HR))$ is given by a flat continuous functor $F= C^\infty(HR)\longrightarrow \Se$. The uniqueness of morphisms in the category $C^\infty(HR)$ entails that a flat functor $F$ to sets only takes two possible values, the empty set and the one element set. Moreover the class $\cN$ of objects $\Omega$ of $C^\infty(HR)$ for which $F(\Omega)\neq \emptyset$ is hereditary and stable under intersection:
$$
\Omega\in \cN, \ \Omega\subset \Omega'\Longrightarrow \Omega'\in \cN, \ \Omega, \, \Omega' \in \cN
\Longrightarrow \Omega\cap\Omega'\in \cN.
$$
Let then $Z(\cN):=\{x\in I\mid \{x\}^c\in \cN\}$. One has $Z(\cN)\neq \emptyset$ if $\cN\neq \{I\}$. Moreover
$$
\Omega\in \cN\iff \Omega^c\subset Z(\cN).
$$
Indeed, if $\Omega\in \cN$ then any element $x$ of its finite complement in $I$ is the complement of a larger $\Omega'\in \cN$ so that $x\in Z(\cN)$. Conversely, if the finite set $\Omega^c\subset Z(\cN)$ then $\Omega$ is a finite intersection of elements of $\cN$ and hence is in $\cN$. Let $E=Z(\cN)^c$. We show that if the flat functor $F$ is continuous then $E$ is convex. Assume, on the contrary, that there exist three elements $x<y<z$ of the interval $I$ such that $x,z\in E$ and $y\notin E$. One has $y\in Z(\cN)$ and thus the open set $\Omega:=\{y\}^c\in \cN$.
\begin{figure}	\begin{center}
\includegraphics[scale=0.6]{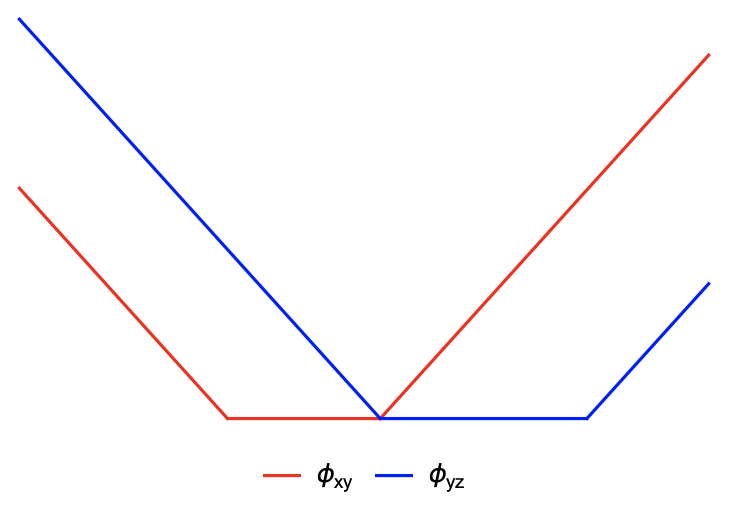}
\end{center}
\caption{$\phi_{xy}\vee \phi_{yz}=\phi_{yy}$ \label{max} }
\end{figure}
 We  construct a covering $(\Omega_j)_{j=1,2}$ of $\Omega$ for the Grothendieck topology $\cJ(HR)$ such that $\Omega_j\notin \cN$ for $j=1,2$ thus contradicting the continuity of the flat functor $F$. \newline
For $a,b\in I$ and $a\leq b$, let $\phi_{ab}$ be the element of $R=C_{\rm Conv}(I)$ which is identically $0$ on the interval $[a,b]$ and has slope $-1$ for $x<a$, and slope $1$ for $x>b$. By construction one has: $Z(\phi_{ab})=\{a,b\}$. Let $x,y,z$ be as above, then (see Figure \ref{max})
$
\phi_{xy}\vee \phi_{yz}=\phi_{yy}
$. This equality shows that the open sets $\Omega_1:=\{x,y\}^c$ and $\Omega_2:=\{y,z\}^c$ form a covering of $\Omega:=\{y\}^c$ for the Grothendieck topology $\cJ(HR)$. 
One has $x\in E=Z(\cN)^c$ and hence $\Omega_1\notin \cN$, similarly $z\in E=Z(\cN)^c$ and hence $\Omega_2\notin \cN$. But since $y\notin E$ one has $\Omega:=\{y\}^c\in \cN$. Thus one obtains that the flat functor $F$ fails to be continuous. This shows that for any continuous flat functor $F$ the subset $E\subset I$ is convex. Conversely, given a convex subset  $E\subset I$ the associated flat functor defined by 
$$
F(\Omega)=\begin{cases} \{*\} &\mbox{if } E\subset \Omega \\ 
\emptyset  & \mbox{otherwise }\end{cases}   
$$
is continuous since $E\subset \Omega$ implies that there exists a connected component $C$ of $\Omega$ which contains $E$. Finally the existence of morphisms of functors is governed by the reverse inclusion of convex sets thus one gets the required result. \endproof \vspace{.05in}

\subsection{The case of semirings}\vspace{.05in}

The tropical case developed in \S\ref{secttropical} together with Proposition 3.4.1 of \cite{SS} suggest that, for an arbitrary semiring $R$, one can compare the spectrum of the $\sss$-algebra $H R$ with
the prime spectrum of $R$.

\begin{prop}\label{semiring1} Let $R$ be a semiring. \newline
$(i)$~The points of the topos $\sspec(HR)$ coincide with the prime spectrum, \ie the set of prime ideals of $R$.	\newline
$(ii)$~The site $(C^\infty(HR),\cJ(HR))$ is the same as the Zariski spectrum $\Spec(R)$ endowed with the basis of open sets of the form $D(f)$.	
\end{prop}
\proof  Once again: a  point of the topos $\sspec(HR)$  is  given by a flat continuous functor $F= C^\infty(HR)\longrightarrow \Se$. The uniqueness of morphisms in the category $C^\infty(HR)$ entails that a flat functor $F$ to sets only takes two possible values: the empty set and the one element set. Moreover, the subset $M\subset R$ of elements $f\in R$ such that  $F(f^\infty)\neq \emptyset$ is hereditary and stable under product:
$$
f\in M\, \&\, \exists n\in \N, \  f^n \in gR \Longrightarrow g\in M, \  \ f,g \in M
\Longrightarrow fg \in M.
$$
Let $J=M^c$ be the complement of $M$. Let $f,g\in J$ and assume that $h=f+g\notin J$.  Lemma \ref{semirings} shows that the pair $\{(fh)^\infty\to h^\infty,(gh)^\infty\to h^\infty\}$ forms a partition of $h^\infty$ (since $h^2=hf+hg$). Thus, if the functor $F$ is continuous, one gets $fh\in M$ or $gh\in M$. Assume $fh\in M$, then as $M$ is hereditary one has $f\in M$, but this contradicts $f\in J$. We have thus shown that $J$ is stable under addition. It is an ideal since $M$ is hereditary and it is a prime ideal since its complement  $M$ is stable under multiplication. Conversely, let $J\subset R$ be a prime ideal. Let $F_J:C^\infty(HR)\longrightarrow \Se$ be the functor defined by 
$$
F_J(f^\infty)=\begin{cases} \{*\} &\mbox{if } f\notin J \\ 
\emptyset  & \mbox{otherwise. }\end{cases}   
$$
Since $J$ is a prime ideal one obtains the flatness of $F_J$. It remains to show that $F_J$ is continuous. By Lemma \ref{semirings}, a family $X$ of morphisms with codomain $f$ is a covering for the  Grothendieck topology $\cJ(A)$ if and only if there exists a map $a:X \to R$ with finite support such that $\sum a(x) x$ is equal to a power of $f$. Assuming that $f\in M$, one needs to show that one has $x\in M$ for some $x\in X$ but this follows since $J$ is an ideal. Thus we have thus shown that the points of the topos $\sspec(HR)$ correspond canonically to prime ideals $\ffp\in {\rm Spec}(R)$ \ie elements of the prime spectrum ${\rm Spec}(R)$. \newline
$(ii)$~We first extend Proposition \ref{hartshorne}  in the semiring context. Proposition 7.28 of \cite{Golan}, states that for any ideal $I\subset R$ of a (commutative) semiring $R$ the intersection of prime ideals containing $I$ is
$$
\sqrt I:=\{ x\in R\mid \exists n\in \N, \ x^n\in I\}.
$$
Setting $V(b):=\{\wp\in {\rm Prim} R\mid b\in \wp\}$ one has 
$$
\sqrt {bR}=\cap_{\wp\in V(b)}\,\wp, \ \sqrt {bR}\supset \sqrt {cR} \iff V(b)\subset V(c), \ \sqrt {bR}\supset \sqrt {cR}\iff \exists n\in \N, \  c^n \in bR
$$
Hence by Proposition \ref{atmostone}  $(iii)$ one gets
$
 V(b)\subset V(c)\iff \exists F_b\to F_c
$.
As in Proposition \ref{hartshorne} it follows that the map $D(b)\mapsto F_b$ is an isomorphism of the category of open sets of ${\rm Spec}( R)$ of the form $D(b)=V(b)^c$, $b\in R$, with the category $C^\infty(HR)$. The open sets $D(b)=V(b)^c$, $b\in R$ form a basis for the Zariski topology and are stable under finite intersections since $D(ab)=D(a)\cap D(b)$. It remains to show that the covering families (say $\{D(b_j)\}$) of an open set $D(b)$) for the  Grothendieck topology $\cJ(HR)$ are the same as for the Zariski topology. By Lemma \ref{semirings}, $\{D(b_j)\}$ is a covering of $D(b)$ for $\cJ(HR)$ if and only if there exist $a_j\in R$ and $n\in \N$  with $b^n=\sum a_jb_j$. If this holds and $\wp\in D(b)$ is a prime ideal, one has $\wp\in D(b_j)$ for some $j$ since otherwise all $b_j\in \wp$ and hence $b^n\in \wp$. Conversely, let $\{D(b_j)\}$ be a covering of $D(b)$ for the Zariski topology. Let $I$ be  the ideal generated by the $b_j$, then any prime ideal $\wp$ containing $I$ is in $\cap V(b_j)\subset V(b)$. Thus the intersection of prime ideals containing $I$ contains $b$ and Proposition 7.28 of \cite{Golan} then shows that there exist $a_j\in R$ and $n\in \N$  with $b^n=\sum a_jb_j$. Thus, the open coverings are the same for the  Grothendieck topology $\cJ(HR)$ as for the Zariski topology.
\endproof

\section{The structure sheaf on $\sspec(A)$}\label{sectstructsheaf}

 In this section we define the structure sheaf on the spectrum of an arbitrary $\sss$-algebra $A$. The construction agrees with the classical one for rings as shown in \S \ref{sectsheafring}. The new feature here is that it is in general false  that this structure presheaf is automatically a sheaf (see Remark \ref{remsheafify}) thus one needs to pass to the associated sheaf: a sheaf of $\sss$-algebras as shown in \S \ref{sectsheafoplusplus}. This construction is based on the process of localization of $\sss$-algebras as explained in \S \ref{sectlocsalg}. The functoriality is proved in \S \ref{sectfunctstrucsheaf}. The new behavior appears  through the operation of quotient by a subgroup and is explored in \S \ref{sectquotG}.  The last \S \ref{sectspech1} investigates the spectrum of the $\sss$-algebra $\vert H\R\vert_1$. \vspace{.05in}
 
 \subsection{Localization of $\sss$-algebras}\label{sectlocsalg}\vspace{.05in}
 
Let $A$ be an $\sss$-algebra and $M\subset A(1_+)$ a multiplicatively closed subset of the monoid $A(1_+)$. Next, we define the localization $M^{-1}A$ as an $\sss$-algebra. 

\begin{lem}\label{quotients0} Let $M$ be a multiplicative, commutative and pointed  monoid. \newline
$(i)$~Let $X$ be an $M$-set, \ie a pointed set endowed with an action of $M$. The following  defines an equivalence relation on $M\times X$ whose quotient is denoted $M^{-1}X$
\begin{equation}\label{defequiv}
	(s,x)\sim (t,y)\iff \exists a\in M\mid atx=asy.
\end{equation}
$(ii)$~Let $f:X\to Y$ be an equivariant map of $M$-sets. Then there exists a unique induced map $\tilde f : M^{-1}X\to M^{-1}Y$ by setting $\tilde f(s,x)=(s,f(x))$ for all $s\in M, x\in X$. \newline
$(iii)$~The association $X\mapsto M^{-1}X$, $f\mapsto \tilde f$ defines a functor $M^{-1}$ from $M$-sets to pointed sets.	
\end{lem}
\proof $(i)$~The relation is symmetric and reflexive by construction. We show that it is transitive. Assume that $(s,x)\sim (t,y)$ and $(t,y)\sim (u,z)$. One then has for some $a,b\in M$, $atx=asy$ and $buy=btz$. It follows, by substitution, that $abtux=abusy=abtsz$, so that $(s,x)\sim (u,z)$.\newline 
$(ii)$~One checks, using the equivariance of $f$, that $(s,x)\sim (t,y)$ implies $(s,f(x))\sim (t,f(y))$.\newline
$(iii)$~follows from $(ii)$.
\endproof

\begin{lem}\label{quotients} Let $A$ be an $\sss$-algebra and $M$ a multiplicatively closed subset of the monoid $A(1_+)$.\newline
$(i)$~The composition $M^{-1}\circ A$ of the functor $M^{-1}$ with $A$ (viewed as a functor  $A:\gop\longrightarrow M$-sets),  defines an $\sss$-module $M^{-1}A$.\newline
$(ii)$~The structure of $\sss$-algebra of $A$ induces a structure of $\sss$-algebra on the $\sss$-module $M^{-1}A$.\newline
$(iii)$~Let $N\subset M$ be a multiplicatively closed subset. Then the inclusion $i:N\subset M$ induces a morphism of $\sss$-algebras: 
$$
N^{-1}A\to M^{-1}A, \ (n,x)\to (i(n),x).
$$
\end{lem}
\proof $(i)$~By \eqref{natu}, the product $A(1_+)\wedge A(k_+)\to A(k_+)$ restricted to $M\wedge A(k_+)$ turns $A$  into a functor  $\gop\longrightarrow M$-sets.  The statement then follows from Lemma \ref{quotients0} $(iii)$.\newline 
$(ii)$~One defines the product $M^{-1}A(X)\wedge M^{-1}A(Y)\to M^{-1}A(X\wedge Y)$ by 
\begin{equation}\label{productloc}
	(s,x)\wedge (t,y)\mapsto (st,m_A(x\wedge y)).
\end{equation}
One checks, using $m_A(ux\wedge vy)= uv \, m_A(x\wedge y)$, $\forall u,v\in M$, that the right hand side only depends upon the classes of $(s,x)$ and $(t,y)$. The naturality, commutativity and associativity of the product follow from the same properties for $A$ and the commutativity of the monoid $M$.\newline
$(iii)$~By construction, the map $(n,x)\mapsto (i(n),x)$ is compatible with the equivalence relation \eqref{defequiv}, since $N\subset M$ and in agreement with the structure of $\sss$-algebras. 
\endproof 

Next proposition shows that the above localization of $\sss$-algebras coincides with the classical localization for semirings and  monoids.

\begin{prop}\label{localok}$(i)$~Let $R$ be a semiring and $M$ a multiplicative subset of $R$. Then the map   $(m,(a_j))\mapsto ((m,a_j))$ defines an isomorphism of $\sss$-algebras $\rho:M^{-1}HR\stackrel{\sim}{\to} HM^{-1}R$.\newline
$(ii)$~Let $N$ be a pointed commutative monoid and $M$ a multiplicative subset of $N$. Then one has $M^{-1}\sss N\simeq \sss M^{-1}N$.
\end{prop}

\proof $(i)$~The map $\rho$ is well defined since for $(m,(a_j))\sim (n,(b_j))$ there exists $u\in M$ such that $(una_j)=(umb_j)$  for all $j$, and this implies $(m,a_j)\sim (n,b_j)$. Conversely, if $(m,a_j)\sim (n,b_j)$  for all $j$, there exists for each $j$ an element $u_j\in M$ such that $u_jna_j=u_jmb_j$, and taking $u=\prod u_j$ one obtains that $(m,(a_j))\sim (n,(b_j))$. This shows that  $\rho$ is injective. It is surjective since for any finite family $(n_j,a_j)$ one has by ``reduction to the same denominator", with $n=\prod n_j$, that for each $j$, $(n_j,a_j)\sim (n,b_j)$, where $b_j=a_j\prod_{i\neq j} n_i$. Finally,  $\rho$ is a morphism of functors $\rho:M^{-1}HR\longrightarrow HM^{-1}R$ and is compatible with the product \eqref{productloc}.\newline
$(ii)$~The statement follows because both $M^{-1}\sss N(k_+)$ and $\sss M^{-1} N(k_+)$ are given by the smash product $M^{-1} N\wedge k_+$. \endproof \vspace{.05in}

\subsection{The sheaf $\cO^{++}$ of $\sss$-algebras}\label{sectsheafoplusplus}\vspace{.05in}
 
 Let $A$ be an $\sss$-algebra. First, we construct a presheaf $\cO$ of $\sss$-algebras for the small category $C^\infty(A)$. Let $f\in A(1_+)$ and $M(f^\infty)$ the multiplicative monoid 
 \begin{equation}\label{structsheaf0}
 	 M(f^\infty):=\{g\in A(1_+) \mid \exists n\in \N, n>0, \ f^n\in gA(1_+)\}.
 \end{equation}
 Using the notations of Lemma \ref{quotients}, we define an $\sss$-algebra   by localization
 \begin{equation}\label{structsheaf}
 \cO(f^\infty):=M(f^\infty)^{-1}A.
 \end{equation}
 We construct the  presheaf $\cO$ as a contravariant functor from the category $C^\infty(A)$ to the category of $\sss$-algebras. 
 Given a morphism $\phi:a^\infty\to b^\infty$ in $C^\infty(A)$, we obtain a morphism of $\sss$-algebras as follows
 
 \begin{lem}\label{structsheaf1} $(i)$~Let $\phi\in \Hom_{C^\infty(A)}(a^\infty,b^\infty)$. The inclusion
 $
  M(b^\infty)\subset  M(a^\infty)
 $
 induces a morphism of $\sss$-algebras $\cO(\phi):\cO(b^\infty)\to \cO(a^\infty)$.\newline
 $(ii)$~The maps  $f^\infty\mapsto \cO(f^\infty)$ and $\phi\mapsto \cO(\phi)$  given in $(i)$, define a contravariant functor from the category $C^\infty(A)$ to the category of $\sss$-algebras.
 \end{lem}
 \proof $(i)$~By Proposition \ref{atmostone} $(iii)$ one has 
 $$
 g\in M(f^\infty)\iff \Hom_{C^\infty(A)}(f^\infty,g^\infty)\neq \emptyset.
 $$
 Thus, the existence of $\phi\in \Hom_{C^\infty(A)}(a^\infty,b^\infty)$ implies $M(b^\infty)\subset  M(a^\infty)$. One then applies Lemma \ref{quotients} $(iii)$ to conclude.\newline  
 $(ii)$~follows from the transitivity of the inclusions $M(c^\infty)\subset M(b^\infty)\subset  M(a^\infty)$.\endproof

 The presheaf $\cO$ constructed in Lemma \ref{structsheaf1} is not always a sheaf as shown in Remark \ref{remsheafify},  thus  one needs to take the associated sheaf. In general, as explained in \cite{MM}, Section III. 5, passing to the associated sheaf requires  two applications of the operation $\cO\mapsto \cO^+$ from presheaves to presheaves (the first  application makes the presheaf  separated and the second makes it a sheaf). As a preliminary step we  show that the operation $\cO\mapsto \cO^+$ transforms presheaves of $\sss$-algebras into presheaves of $\sss$-algebras. This amounts to prove  that the matching families for a given sieve still form an $\sss$-algebra. Next we recall the definition of matching families (\!\cite{MM}, p. 121).
 
 \begin{defn}\label{defnmatchfam}
 	Given a site, a sieve $R$ which is a cover of an object $C$ and a presheaf $P$, a matching family is a function which assigns to every element $f:D\to C$ of $R$ an element $x_f\in P(D)$ so that 
\begin{equation}\label{matching0}
 P(g)(x_f)=x_{fg}\qqq g:E\to D.	
\end{equation}
 \end{defn} 
 
 One denotes by ${\rm Match}(R,P)$ the matching families for the sieve $R$ and the presheaf $P$.\vspace{.05in}
 
 Let $(A_\alpha)_{\alpha \in I}$ be a family of $\sss$-algebras: we define the product $\sss$-module $A:=\prod A_\alpha$  as 
 $$
 A(k_+):=\prod A_\alpha(k_+), \ A(\phi)((a_\alpha)_{\alpha \in I}):=(A_\alpha(\phi)(a_\alpha))_{\alpha \in I}.
 $$
  The product $m_A:A(X)\wedge A(Y)\to A(X\wedge Y) $ is given by 
 $$
 m_A((a_\alpha)_{\alpha \in I},(b_\alpha)_{\alpha \in I}):=(m_{A_\alpha}(a_\alpha,b_\alpha))_{\alpha \in I}.
 $$
 This endows $A$ with a structure of $\sss$-algebra. 
 \begin{lem}\label{matching} Let $\cC$ be a small category, $C$ an object of $\cC$, $R$ a sieve on $C$ and $P$ a presheaf of $\sss$-algebras. Let $P_k$ be the presheaf of sets obtained by evaluation on the object $k_+$ of $\gop$. Then the map $k_+\to {\rm Match}(R,P_k)$ defines an $\sss$-subalgebra of the product $\sss$-algebra: 
 $A=\prod_{f\in R} P({\rm Dom}(f))
 $.  	
 \end{lem}
\proof Let us first show that $k_+\to {\rm Match}(R,P_k)$ is a subfunctor of the functor $k_+\to A(k_+)$. By construction, the matching families $(x_f)_{f\in R}\in {\rm Match}(R,P_k)$ form a subset of $A(k_+)$. Let $\phi\in \Hom_\gop(k_+,\ell_+)$ and $(x_f)_{f\in R}\in {\rm Match}(R,P_k)$. We show that $A(\phi)(x_f)_{f\in R}\in {\rm Match}(R,P_\ell)$. Let $g:E\to D$ be a morphism in $\cC$. Since $P(g):P(D)\to P(E)$ is a morphism of $\sss$-algebras, it is a natural transformation of functors $\gop \longrightarrow \Se$ and this gives the commutative diagram 
 \begin{gather} \raisetag{37pt} \,\hspace{60pt}
\xymatrix@C=25pt@R=25pt{
P_k(D) \ar[d]_{P_k(g)} \ar[r]^{P(D)(\phi)} &
P_\ell(D)\ar[d]^{P_\ell(g)} \\
 P_k(E) \ar[r]^{P(E)(\phi)}  & P_\ell(E)
  } \label{natuC}\hspace{100pt} 
\end{gather}
 One has $x_f\in P_k(D)$,  $A(\phi)(x_f)=P(D)(\phi)(x_f)$, and thus, using \eqref{natuC} and \eqref{matching0} one obtains
 $$
 P_\ell(g) A(\phi)(x_f)= P_\ell(g)P(D)(\phi)(x_f)=P(E)(\phi) P_k(g)x_f=P(E)(\phi) x_{fg}=A(\phi)(x_{fg}).
 $$
 Let us show that $k_+\to {\rm Match}(R,P_k)$ is stable under the  product $m_A:A(k_+)\wedge A(\ell_+)\to A(k_+\wedge \ell_+)$. Since $P(g):P(D)\to P(E)$ is a morphism of $\sss$-algebras one has a commutative diagram 
 \begin{gather} \raisetag{37pt} \,\hspace{60pt}
\xymatrix@C=25pt@R=25pt{
P_k(D)\wedge P_\ell(D) \ar[d]_{P_k(g)\wedge P_\ell(g)} \ar[r]^{\ \ m_A} &
P_{k\ell}(D) \ar[d]^{P_{k\ell}(g)} \\
P_k(E)\wedge P_\ell(E) \ar[r]_{\ \ m_A}  & P_{k\ell}(E)
  } \label{diag2ter}\hspace{100pt}
\end{gather}
Let then $(x_f)_{f\in R}\in {\rm Match}(R,P_k)$ and $(y_f)_{f\in R}\in {\rm Match}(R,P_\ell)$. To show that $m_A(x_f\wedge y_f)_{f\in R}\in {\rm Match}(R,P_{k\ell})$ one uses
$$
P_{k\ell}(g)(m_A(x_f\wedge y_f))=m_A(P_k(g)x_f\wedge P_\ell(g)y_f)=m_A(x_{fg}\wedge y_{fg}).
$$
This proves that the map $k_+\to {\rm Match}(R,P_k)$ defines an $\sss$-subalgebra of $A$.\endproof 

 \begin{lem}\label{matching2} Let $(\cC,J)$ be a Grothendieck site,  and $P$ a presheaf of $\sss$-algebras. Then the associated sheaf $k_+\to P_k^{++}$ is a sheaf of $\sss$-algebras.
 \end{lem} 
 \proof The associated sheaf is obtained by iterating twice the transformation $P\mapsto P^+$ on presheaves, thus it suffices to show that this transformation maps presheaves of $\sss$-algebras to presheaves of $\sss$-algebras. One has by construction
 $$
  P_k^+(C):=\varinjlim_{R\in J(C)}{\rm Match}(R,P_k)
 $$
 so that, by Lemma \ref{matching}, $P^+(C)$ is a pointwise  colimit of $\sss$-algebras and is thus an $\sss$-algebra. \endproof 
 
  \begin{defn}\label{specAbis} Let $A$ be an $\sss$-algebra. The {\em structure presheaf} $\cO_A$ (resp. {\em structure  sheaf} $\cO_A^{++}$) on the site $(C^\infty(A),\cJ(A))$  is the  presheaf $\cO$ of Lemma \ref{structsheaf1} (resp. the sheaf $\cO^{++}$ of $\sss$-algebras associated to $\cO$ on $\sspec(A)$).
 \end{defn}\vspace{.05in}

\subsection{The case of rings}\label{sectsheafring}\vspace{.05in}

 In this section we prove that the general construction of the structure sheaf on $\sspec(A)$ as presented in \S\ref{sectsheafoplusplus} agrees with the classical one when $A=HR$ for a ring $R$. 
 
 \begin{prop}\label{sheafR} Let $R$ be a ring and $A=HR$ be the associated $\sss$-algebra. Then the structure presheaf  $\cO_A$ on  $\sspec(A)$ is a sheaf and is  canonically isomorphic to the sheaf of $\sss$-algebras $H\cO_R$ on $\Spec R$. 	
 \end{prop}
 \proof It is well known (\!\cite{Hart}, Proposition II. 2.2) that for an element $f\in R$ the ring $\cO_R(D(f))$, for the classical structure sheaf $\cO_R$ on the prime spectrum, is isomorphic to the localized ring $R_f$. This latter ring is the same as the localization $M(f^\infty)^{-1}R$, since the elements involved in $M(f^\infty)$ in \eqref{structsheaf0} all divide a power of $f$. By Proposition \ref{localok}, one derives  $HR_f=M(f^\infty)^{-1}HR$ and one has by \opcit  that the presheaf $D(f)\mapsto R_f$ is already a sheaf. It follows that the presheaf on $\sspec(A)$ given by $f^\infty \mapsto HR_f$ is a sheaf, and thus so is the presheaf $f^\infty \mapsto M(f^\infty)^{-1}HR$.  This shows that this presheaf coincides with $\cO_A$ and that, in the above case, there is no need to sheafify the presheaf of Lemma \ref{structsheaf1}.\endproof\vspace{.05in} 

\subsection{Functoriality of the  structure sheaf $\cO_A^{++}$}\label{sectfunctstrucsheaf}\vspace{.05in}

Let  $\phi: (C,J)\to (D,K)$ be a  morphism of sites as in \S \ref{sectfunctorial}. Let $f:{\rm Sh}(D,K)\longrightarrow {\rm Sh}(C,J)$ be the associated geometric morphism and $f_*$ the direct image functor. We denote by 
 $\mathfrak a$ the sheafification functor.  Let $\cF$ be a presheaf on $D$, then one has a natural morphism of sheaves 
\begin{equation}\label{sheafify}
\gamma: {\mathfrak a}f_*\cF \to f_*{\mathfrak a}\cF
 \end{equation}
induced by the universal property of ${\mathfrak a}f_*\cF$. Indeed,  any map $\psi$ of presheaves from $f_*\cF$ to a sheaf $\cL$ factors uniquely through the canonical map $\eta: f_*\cF\to {\mathfrak a}f_*\cF$ (see \cite{MM} Lemma III. 5.3), so that one gets a morphism of sheaves $\tilde \psi:{\mathfrak a}f_*\cF\to \cL$ such that $\psi=\tilde \psi \circ \eta$.  Here, one has a canonical morphism of presheaves  $\cF\to{\mathfrak a}\cF$ and by functoriality of $f_*$ at the level of presheaves (since it is given by composition with $\phi:C\longrightarrow D$) one gets a corresponding morphism $\psi :f_*\cF\to f_*{\mathfrak a}\cF$. Then one gets $\gamma=\tilde \psi$. 

Let now $\cE$ be a presheaf on $C$ and $\rho:\cE\to  f_*\cF$ a morphism of presheaves. By functoriality of the sheafification functor $\mathfrak a$ one obtains a morphism of sheaves $\mathfrak a(\rho):\mathfrak a(\cE)\to \mathfrak a(f_*\cF)$. We let
\begin{equation}\label{sheafify1}
\alpha( \rho):=\gamma \circ \mathfrak a(\rho):\mathfrak  a(\cE)\to f_*{\mathfrak a}(\cF)
 \end{equation}
be the morphism of sheaves  obtained by composition with $\gamma$ of \eqref{sheafify}. 
\begin{thm}\label{thmfunctorsheaf} Let $A,B$ be $\sss$-algebras and $\phi:A\to B$ a morphism of $\sss$-algebras. Let $\tilde \phi: (C^\infty(A),\cJ(A))\longrightarrow  (C^\infty(B),\cJ(B))$ be the associated morphism of sites as in Theorem \ref{functori}. The localization of $\phi$ defines a morphisms of presheaves $\rho:\cO_A\to \tilde \phi_*(\cO_B)$
\begin{equation}\label{sheafify2}
\rho:M(f^\infty)^{-1}A \to M(\phi(f)^\infty)^{-1}B.
 \end{equation}
Moreover, the  associated morphism of sheaves $\alpha(\rho)$ as in \eqref{sheafify}, is a morphism from the structure sheaf $\cO_A^{++}$ to the direct image of the structure sheaf  $\cO_B^{++}$ by the geometric morphism associated to the morphism of sites $\tilde \phi$.
\end{thm}
 \proof 
 The map which to a pair $(m,a)\in M(f^\infty)\times A$ assigns the pair $(\phi(m),\phi(a))$ that belongs to $ M(\phi(f)^\infty)\times B$, defines, as in Lemma  \ref{quotients} $(iii)$, a morphism of $\sss$-algebras. Moreover these morphisms are compatible with the restriction morphisms of Lemma \ref{structsheaf1}, thus they define a morphism of presheaves as in \eqref{sheafify2}, and $\rho$ is  a morphism of presheaves of $\sss$-algebras. By Lemma \ref{matching} the sheafification functor as well as the universal map from a presheaf to the associated sheaf are compatible with the $\sss$-algebra structures. It follows that the morphism of sheaves $\alpha(\rho)$ is itself compatible with the $\sss$-algebra structures.\endproof \vspace{.05in}

\subsection{Quotient by a subgroup}\label{sectquotG}\vspace{.05in}

Let $A$ be an $\sss$-algebra and $G\subset A(1_+)$ be a multiplicative subgroup of the monoid $A(1_+)$. By \eqref{covariance} the functor $A:\gop\longrightarrow \Se_*$ is $G$-equivariant and induces, by composition with the quotient by $G$, a functor  $A/G:\gop\longrightarrow \Se_*$ which assigns to an object $X$ of $\gop$ the set of $G$-orbits $A(X)/G$ and to a morphism its action on the $G$-orbits. Moreover, the product $m_A:A(X)\wedge A(Y)\to A(X\wedge Y)$ fulfills $m_A(u\xi\wedge v\eta)=uv\ m_A(\xi\wedge \eta)$ and hence induces a product which turns $A/G$ into an $\sss$-algebra. 
We let $q_G:A\to A/G$ be the natural transformation given by the canonical map to the quotient by $G$. It is by construction a morphism of $\sss$-algebras. 

\begin{prop}\label{propquot} Let $A$ be an $\sss$-algebra, $G\subset A(1_+)$ a subgroup and $q_G:A\to A/G$ the canonical morphism to the quotient by $G$. \newline
$(i)$~The geometric morphism $\sspec(A/G)\to \sspec(A)$ induced by $q_G$ is an isomorphism of sites.\newline
$(ii)$~The morphism  $\rho$ of presheaves as in \eqref{sheafify2} is given by the localization of the morphism $q_G$.	
\end{prop}
\proof $(i)$~We show that the morphism of sites $\tilde q_G: (C^\infty(A),\cJ(A))\longrightarrow  (C^\infty(A/G),\cJ(A/G))$ is an isomorphism. First, the morphism $q_G:A\to A/G$ induces a morphism of monoids $q_G:A(1_+)\to A/G(1_+)$ whose effect on $M=A(1_+)$ is to divide $M$ by the subgroup $G$. We claim that the quotient map $q:M\to M/G$ induces an isomorphism $\tilde q_G:C^\infty(A)\stackrel{\sim}{\to} C^\infty(A/G)$. To prove this it is enough to show that for any $f\in M$ and $g\in G$ one has the equality $f^\infty=(fg)^\infty$ in $C^\infty(A)$. This follows from Proposition \ref{atmostone}, since $f$ divides $fg$ and $fg$ divides $f$.  Let us show that $\tilde q_G:C^\infty(A)\to C^\infty(A/G)$ is an isomorphism for the Grothendieck topologies. By Definition \ref{partit}, a partition of an object $f$ in $C^\infty(A)$ is a collection of morphisms $f_j\to f$, $j=1, \ldots ,n$,  such that there exists $\xi\in A(n_+)$ with   $(A(\delta_j)(\xi))^\infty=f_j$, $\forall j$, and $(A(\Sigma) \xi)^\infty=f$. This notion is unchanged if one replaces $A$ by $A/G$. It follows that the notion of multi-partition is also unchanged and the morphism of sites $\tilde q_G: (C^\infty(A),\cJ(A))\longrightarrow  (C^\infty(A/G),\cJ(A/G))$ is an isomorphism.\newline
$(ii)$~ follows from the commutation of the operations of localization with respect to a multiplicative subset $M\subset A(1_+)$ containing $G$, and of quotient by $G$, \ie the equality 
\begin{equation}\label{sheafify3}
(M/G)^{-1}(A/G)=(M^{-1}A)/G.
\end{equation} \endproof 
Next result displays the role of the general theory of $\sss$-algebras in the context of the adele class space (compare $(i)$ with example \S 3.5, example 11, in \cite{compositio}).
\begin{cor}\label{adeleclass} Let $K$ be a global field, $\A_K$ the ring of adeles of $K$. Let $A=H\A_K/K^\times$ be the $\sss$-algebra associated as above to the subgroup $G=K^\times\subset H\A_K(1_+)=\A_K$. \newline
$(i)$~The map of sets $\nu:\Hom_\salg(\sss[T],A)\to A(1_+)$ of Corollary \ref{transit} defines a canonical bijection 
$$
\nu:\Hom_\salg(\sss[T],H\A_K/K^\times)\stackrel{\sim}{\to} \A_K/ K^\times.
$$
$(ii)$~The spectrum $\sspec(H\A_K/K^\times)$ is canonically isomorphic as a site with the spectrum of the ring of adeles  	$\A_K$.
\end{cor}
\proof $(i)$~follows directly from Corollary \ref{transit} and the equality $H\A_K(1_+)=\A_K$.\newline
$(ii)$~follows from Proposition \ref{propquot} and Proposition \ref{sheafR}. \endproof

\begin{rem}\label{remsheafify} 
It is false, in general and in the above context, that the presheaf given by \eqref{sheafify3} on $C^\infty(A/G)$ is a sheaf. This shows, in particular, that one needs to pass to the associated sheaf. In the following we describe some examples where this issue happens. \newline
Let $R=R_1\times R_2$ be the product of two commutative rings, $G$ a commutative group and $\iota_j:G\to R_j$ a morphism to the group of units, for $j=1,2$. To be more specific, we consider two places $v_j$ of $\Q$, we let $R_j=\Q_{v_j}$ and $G=\Q^\times$ with the canonical inclusions in the local fields obtained as completions (much of this argument applies in a more general setup). The spectrum of $HR/G$ is, by Propositions \ref{ordinary} and \ref{propquot} the same as the prime spectrum of $R$ and hence consists of the two point set $\{v_1,v_2\}$ endowed with the discrete topology. The structure sheaf of $HR$ is obtained by applying $H$ to the structure sheaf of $R$ and the latter has stalk $R_j$ at $v_j$. The canonical presheaf on $\sspec(HR/G)$ assigns to  the open set $\{v_1,v_2\}$ the $\sss$-algebra $HR/G$, but when one passes to the associated sheaf one obtains a non-trivial quotient. One can see this at the level $1_+$, \ie by evaluation on $1_+$.  For the canonical presheaf $\cO$ one obtains
$$
\cO(\{v_1,v_2\})(1_+)=(R_1\times R_2)/G, \ \cO(\{v_1\})(1_+)=R_1/G, \ \cO(\{v_2\})(1_+)= R_2/G,
$$
so that for the associated sheaf one has 
$$
\cO_A(\{v_1,v_2\})(1_+)=R_1/G\times R_2/G.
$$
This shows that the canonical morphism to the associated sheaf is in general not injective. 
	
\end{rem}
In view of the above Remark \ref{remsheafify}  the canonical presheaf on  $\sspec(HR/G)$ is in general not a sheaf. However, next statement guarantees that it is a sheaf provided the ring $R$ has no zero divisors.

\begin{prop}\label{propquot1}
Let $R$ be a commutative ring with no zero divisors, and $G\subset R$ a subgroup of the multiplicative group. Then the 	canonical presheaf on  $\sspec(HR/G)$ is  a sheaf.
\end{prop}
\proof By Propositions \ref{ordinary} and \ref{propquot}, $\sspec(HR/G)$ is the same as the prime spectrum of $R$. It is thus an irreducible topological space, \ie such that any pair of non-empty open sets has a non-empty intersection. Let $\cO_R$ be the structure sheaf of $R$. Since $R$ has no zero divisors the restriction maps $\cO_R(U)\to \cO_R(V)$ are injective and the action of $G$ by multiplication is free except for its action on the fixed point $0$. These two properties continue to hold for finite powers $\cO_R^n$ and are hence fulfilled by the structure sheaf of $HR$. Then the result is ensured by the next lemma.\endproof 

\begin{lem}\label{h1irred} Let $X$ be an irreducible topological space, $G$ a group and $\cF$ a sheaf of pointed $G$-sets on $X$ such that 
\begin{enumerate}
\item The action of $G$ on the complement of the base point is free on $\cF(U)$ for any non-empty open set $U$.
\item The restriction maps 	$\cF(U)\to \cF(V)$ are injective.
\end{enumerate}
Then the presheaf $U\mapsto \cF(U)/G$ is a sheaf.	
\end{lem}

\proof Let $(U_j)_{j\in I}$ be an open cover of the open set $U\subset X$ and consider the equalizer diagram for the sheaf $\cF$
\begin{equation}\label{equalizer}
	\cF(U)\stackrel{\epsilon}{\to} \prod_{j\in I} \cF(U_j)\stackrel{\rightarrow }{\rightarrow}\prod_{i,j\in I} \cF(U_i\cap U_j).
\end{equation}
Let $i_0\in I$, with $U_{i_0}\neq \emptyset$.  The restriction map $\cF(U)\to \cF(U_{i_0})$ is injective by $(ii)$, and $G$-equivariant, thus it remains injective on the orbit space \ie for the induced map 
$\cF(U)/G\to \cF(U_{i_0})/G$. Let $\xi=(\xi_j)_{j\in I}\in \prod_{j\in I} \cF(U_j)$ be such that the restrictions of $\xi_i$ and $\xi_j$ are equal in $\cF(U_i\cap U_j)/G$ for all $i,j\in I$. We can assume that all $\xi_j$ are distinct from the base point $*$, since if one of them is $*$ the same holds for all. Since the action of $G$ is free on  the complement of the base point there exists unique elements $g(i,j)\in G$ such that 
$$
g(i,j)\xi_j=\xi_i \ \text{on} \ U_i\cap U_j\qqq i,j\in I. 
$$
One has $g(i,i)=1$, $g(i,j)g(j,k)=g(i,k)$ for all $i,j,k\in I$. Thus $g(i,j)=g(i,i_0)g(j,i_0)^{-1}$ and with $\eta_j:=g(j,i_0)^{-1}\xi_j$ one has 
$$
\eta_j=\eta_i \ \text{on} \ U_i\cap U_j\qqq i,j\in I,
$$
so that since \eqref{equalizer} is an equalizer there exists a section $\eta\in \cF(U)$ which restricts to $\eta_j$ on $U_j$ for all $j\in I$. This shows that the class of $\eta$ in $\cF(U)/G$ restricts to the class of $\xi_j$ in $\cF(U_j)/G$ for all $j\in I$. Hence the diagram \eqref{equalizer} remains an equalizer after passing to the sets of $G$-orbits. \endproof \vspace{.05in}

\subsection{The spectrum of $\vert H\Q\vert_1$}\label{sectspech1}\vspace{.05in}

As in \cite{CCprel}, we consider the $\sss$-subalgebra of $H\Q$ defined by 
$$
\vert H\Q\vert_1(n_+):=\{ (q_j)_{j=1,\ldots,n}\in H\Q(n_+)\mid \sum \vert q_j\vert\leq 1\},
$$
where $\vert q\vert$ is the archimedean norm of $q\in \Q$. Next proposition shows that $\sspec(\vert H\Q\vert_1)$ and its structure sheaf behave in a similar way as for the local rings associated to the non-archimedean places of $\Q$, but with a substantial nuance for the Grothendieck topology.

\begin{prop}\label{prophq1}$(i)$~The spectrum $\sspec(\vert H\Q\vert_1)$ is the site given by  the small category $0\to u\to 1$ endowed with the Grothendieck topology for which the single morphism $u\to 1$ is a  covering of $1$. \newline
$(ii)$~The   structure presheaf $\cO$ is such that  $\cO(1)=\vert H\Q\vert_1$, $\cO(u)= H\Q$  and the restriction map is given by the inclusion $\vert H\Q\vert_1\subset H\Q$.	
\end{prop}
\proof $(i)$~Let $A=\vert H\Q\vert_1$. One has $A(1_+)=\{ q\in \Q\mid  \vert q\vert\leq 1\}$. Thus for $r,s\in A(1_+)$, one has
$$
r\in s A(1_+)\iff \vert r\vert \leq \vert s\vert.
$$
This proves that the elements $u=r^\infty$ for $0<\vert r\vert<1$ are all equal in the  category $C^\infty(A)$. This category has  three objects, $u$, $0^\infty$ and $1^\infty$. The morphisms are $0^\infty\to u \to 1^\infty$. The partition of unity given by $\frac 12 +\frac 12=1$ shows that the single morphism $u\to 1$ is a  covering of $1$. \newline
$(ii)$~The localization of the $\sss$-algebra $A=\vert H\Q\vert_1$, with respect to the monoid $M$ of non-zero elements of $A(1_+)$, is $H\Q$ since for any $n$-tuple $(q_j)_{j=1,\ldots,n}\in H\Q(n_+)$ one can find a non-zero multiple $(qq_j)_{j=1,\ldots,n}\in H\Q(n_+)$ such that $\sum \vert qq_j\vert\leq 1$.\endproof 

 In the above proposition $u$ corresponds to the generic point and the stalk at $u$ is $H\Q$ as expected, but the closed point disappears due to the different Grothendieck topology. This example exhibits the need to keep all the information from the site and presheaf instead of passing directly to the associated topos and sheaf.

\section{Relation with T\"oen-Vaqui\' e \cite{TV}}\label{sectfinalrems}

In their paper \cite{TV}, B. T\"oen and M.  Vaqui\' e have developed a general theory of algebraic geometry which applies to any symmetric monoidal closed category $C$ which is complete and cocomplete. Their theory gives the expected result, \ie usual algebraic geometry, when applied to the category of abelian groups. It also gives the naive $\F_1$ theory of monoids when applied to the category of pointed sets as suggested by \cite{KS}. Now, the category $\smod$ of $\Gamma$-sets (equivalently of $\sss$-modules) is a symmetric monoidal closed category  which is complete and cocomplete.  The closed structure of the category $\smod$  of $\Gamma$-sets  is defined by setting
\begin{equation}\label{closedstructure0}
	\underline\Hom_\sss(M,N)=\{k_+\mapsto \Hom_\sss(M,N(k_+\wedge -))\},
\end{equation}
where $\wedge$ is the smash product of pointed sets. This formula uniquely defines the smash product of $\Gamma$-sets by applying the adjunction 
$$
\underline\Hom_\sss(M_1\wedge M_2,N)=\underline\Hom_\sss(M_1,\underline\Hom_\sss(M_2,N)).
$$
Thus it would seem  perfectly natural to apply the general theory of \cite{TV} to the category $\smod$ rather than developing the new theory presented above. However we shall prove now that even though rings  form a full subcategory   $\An$ of the category of $\sss$-algebras, the theory as in \cite{TV}, when applied to the category $\smod$,  does not agree with ordinary algebraic geometry when restricted to $\An$. 

\begin{lem}\label{twofields} Let $K$ be a field $K\neq \F_2$, $R$ the ring $R=K^2$, and $A=HR$.  Let $p_j:R\to K$, $j=1,2$ be the two projections. Then the family with two elements
$$
(Hp_j)_{j=1,2}, \ \ Hp_j: HR\to HK
$$
fails to be a Zariski covering of $HR$ in the sense of 	\cite{TV}, Definition 2.10. 
\end{lem}
\proof Let $G=K^*$ be the multiplicative group of $K$, viewed as a diagonal subgroup of $R=K^2$, using the morphism $g\mapsto (g,g)\in  K^2$. Let $A':=HR/G$ and $\rho:A\to A'$ the canonical morphism. One has $HR=HK\times HK$ and the morphism $Hp_j: HR\to HK$ is simply the canonical projection to one of the factors. We view $A'$ as an $A$-module using the morphism $\rho$, and we first compute the tensor product $X_j:=A'\otimes_A B_j$, where $A$ acts on $B_j=HK$ using the morphism $Hp_j$. By  \cite{TV} page 454, one has a  natural isomorphism of $B_j$-modules of the form 
$$
A'\coprod_A B_j\to A'\otimes_A B_j,
$$ 
where the colimit $A'\coprod_A B_j$ is computed in the category $\salg$ of $\sss$-algebras. We show that $A'\coprod_A B_j$ is given by 
\begin{gather} \raisetag{37pt} \,\hspace{60pt}
\xymatrix@C=25pt@R=25pt{
& HR \ar[dl]_{\rho} \ar[dr]^{Hp_j} &\\
HR/G \ar[dr]^{\pi_j}& & HK \ar[dl]_{\rho_j} \\
 & HK/K^* &
  } \label{directsum1}\hspace{100pt} 
\end{gather}
 where $\rho_j:HK\to HK/K^*$ is the quotient morphism and $\pi_j:HR/G\to HK/K^*$ the projection. Let $X$ be an object of $\salg$  and morphisms $\alpha:A'\to X$ and $\beta:B_j\to X$ such that the diagram \eqref{directsum} commutes
 \begin{gather} \raisetag{37pt} \,\hspace{60pt}
\xymatrix@C=25pt@R=25pt{
& A \ar[dl]_{\rho} \ar[dr]^{Hp_j} &\\
A' \ar[dr]^\alpha& & B_j \ar[dl]_\beta \\
 & X &
  } \label{directsum}\hspace{100pt} 
\end{gather}
We prove that, for fixed $j$, there exists a  unique morphism $\gamma:HK/K^*\to X$ such that 
$$
\alpha=\gamma \circ \pi_j, \ \ \beta =\gamma \circ \rho_j.
$$
We take $j=1$ for simplicity. By construction $\beta$ is a morphism of $\sss$-algebras from $HK$ to $X$. 
For any integer $n>0$ and any $\xi_1,\xi_2 \in HK(n_+)$ one has 
$(\xi_1,\xi_2)\in  HR(n_+) $, and
$$
\alpha(\rho((\xi_1,\xi_2)))=\alpha((\xi_1,\xi_2)G)=\beta(\xi_1).
$$
This shows that $\beta(\xi_1)$ is unchanged if one multiplies $\xi_1$ by any $u\in K^*$ and that there exists a  unique morphism $\gamma:HK/K^*\to X$ such that $\beta =\gamma \circ \rho_1$. Then one gets $\alpha=\gamma \circ \pi_1$. This shows that 
the colimit $A'\coprod_A B_j$ is given by the diagram \eqref{directsum1}. 
Thus the functor $\cL$
$$
Y\to \prod Y\otimes_A B_j, \ A-{\rm Mod} \longrightarrow \prod B_j-{\rm Mod}
$$
when evaluated on the $A$-module $A'$, transforms $A'$ into the product of two copies of $HK/K^*$, viewed as modules over $B_j=HK$. Let $u\in K^*$, $u\neq 1$ and $v=(1,u)\in R=A(1_+)$. The endomorphism $V$ of  the $A$-module $A'$ given by multiplication by $\rho(v)\in A'(1_+)$ is well defined and non trivial since $\rho(v)\neq 1$. But the action of $V$ on each $A'\otimes_A B_j=HK/K^*$ is the identity and thus the functor $\cL$ is not faithful. This shows that the covering is not faithfully flat and hence is not a Zariski cover. \endproof 

Lemma \ref{twofields} shows that T\"oen-Vaqui\' e
's theory for the   symmetric monoidal closed category $\smod$ cannot restrict to the usual algebraic geometry on the full subcategory $\An$ of $\salg$. Indeed, if it would, the spectrum of a field $K_j$ would be a single point $\{p_j\}$ and the spectrum of the product $R=K_1\times K_2$ of two fields would be a pair of points $\{p_1,p_2\}$ with maps corresponding to the two projections $p_j:R\to K_j$ forming a Zariski cover. On the other hand, Lemma \ref{twofields} shows that one does not obtain a Zariski cover in the sense of the T\"oen-Vaqui\' e theory. The conceptual reason for this failure is that the notion of cover as in \cite{TV} involves all $A$-modules.  Even though ordinary $R$-modules, for a ring $R$, form a full subcategory of the category of $HR$-modules, there are new $HR$-modules in the wider context of $\sss$-algebras: their existence constrains the notion of cover and excludes the simplest covers. 

\section{Appendix}\label{sectadjfunctor} \vspace{.05in}

\subsection{Gluing two categories using adjoint functors}\vspace{.05in}

 We briefly recall, for convenience of the reader, the process  of glueing together two categories using a pair of adjoint functors, as explained to us by P. Cartier (see  \cite{compositio} for more details).\newline
 Let $\cC$ and $\cC'$ be two categories connected by  a pair of adjoint functors $\beta:\cC\to \cC'$
and $\beta^*:\cC'\to \cC$. By definition one has a canonical identification
\begin{equation}\label{adjtrel0}
\Hom_{\cC'}(\beta(H), R)\stackrel{\Phi}{\cong} \Hom_{\cC}(H,\beta^*(R))\qquad\forall~H\in{\rm Obj}(\cC),~R\in{\rm Obj}(\cC').
\end{equation}
The naturality of  $\Phi$ is expressed by the commutativity of the following diagram where the vertical arrows are given by composition, $\forall f\in \Hom_{\cC}(G,H)$ and $\forall h\in \Hom_{\cC'}(R,S)$
\begin{gather} \raisetag{37pt} \,\hspace{60pt}
\xymatrix@C=25pt@R=25pt{
\Hom_{\cC'}(\beta(H),R) \ar[d]_{\Hom(\beta(f),h)} \ar[r]^{\Phi} &
\Hom_{\cC}(H,\beta^*(R)) \ar[d]^{\Hom(f,\beta^*(h))} \\
  \Hom_{\cC'}(\beta(G),S) \ar[r]_{\Phi}  & \Hom_{\cC}(G,\beta^*(S))
  } \label{diag}\hspace{100pt}
\end{gather}

We  define the category $\cC''=\cC\cup_{\beta,\beta^*} \cC'$ by gluing $\cC$ and $\cC'$. The collection of objects of $\cC''$ is obtained as the {\em disjoint union} of the collection of objects of $\cC$ and $\cC'$. For $R\in{\rm Obj}(\cC')$ and $H\in{\rm Obj}(\cC)$, one sets $\Hom_{\cC''}(R,H)=\emptyset$. On the other hand, one defines
\begin{equation}\label{morC0}
    \Hom_{\cC''}(H,R)=\Hom_{\cC'}(\beta(H), R)\cong \Hom_{\cC}(H,\beta^*(R)).
\end{equation}
 The morphisms between objects contained in a same category are unchanged. The composition of morphisms in $\cC''$ is defined as follows. For $\phi \in \Hom_{\cC''}(H,R)$ and $\psi \in \Hom_{\cC}(H',H)$, one  defines $\phi\circ \psi\in \Hom_{\cC''}(H',R)$ as the composite
 \begin{equation}\label{defcomp1}
 \phi\circ \beta(\psi)\in \Hom_{\cC'}(\beta(H'), R)=\Hom_{\cC''}(H',R).
 \end{equation}
  Similarly, for $\theta\in \Hom_{\cC'}(R, R')$ one defines $\theta \circ \phi\in \Hom_{\cC''}(H,R')$ as the composite
  \begin{equation}\label{defcomp2}
 \theta \circ \phi\in \Hom_{\cC'}(\beta(H), R')=\Hom_{\cC''}(H,R').
 \end{equation}
 We recall from \cite{compositio} the following fact
 \begin{prop}\label{catplus0}  $\cC''=\cC\cup_{\beta,\beta^*} \cC'$ is a category which contains $\cC$ and $\cC'$ as full subcategories. Moreover, for any object $H$ of $\cC$ and $R$ of $\cC'$, one has
$$
 \Hom_{\cC''}(H,R)=\Hom_{\cC'}(\beta(H), R)\cong \Hom_{\cC}(H,\beta^*(R)).
$$
\end{prop}
One defines specific morphisms $\alpha_H$ and $\alpha'_R$ as follows
\begin{equation}\label{alphas}
    \alpha_H={\rm id}_{\beta(H)}\in \Hom_{\cC'}(\beta(H), \beta(H))=\Hom_{\cC''}(H,\beta(H))
\end{equation}
\begin{equation}\label{alphasbis}
    \alpha'_R=\Phi^{-1}({\rm id}_{\beta^*(R)})\in \Phi^{-1}( \Hom_{\cC}(\beta^*(R), \beta^*(R)))=\Hom_{\cC''}(\beta^*(R),R).
\end{equation}
By construction one gets
\begin{equation}\label{expressbis}
    \Hom_{\cC''}(H,R)=\{g\circ\alpha_H\,|\, g\in \Hom_{\cC'}(\beta(H),R)\},
\end{equation}
and for any morphism $\rho\in \Hom_{\cC}(H,K)$ the following equation holds
\begin{equation}\label{exprrel}
\alpha_K\circ \rho=\beta(\rho)\circ \alpha_H.
\end{equation}
Similarly, it also turns out that
\begin{equation}\label{express}
    \Hom_{\cC''}(H,R)=\{\alpha'_R\circ f\,|\, f\in \Hom_{\cC}(H,\beta^*(R))\}
\end{equation}
and the associated equalities hold
\begin{equation}\label{exprrelbis}
\alpha'_S\circ \beta^*(\rho)= \rho\circ \alpha'_R \qquad\forall \rho\in \Hom_{\cC'}(R,S)
\end{equation}
\begin{equation}\label{relate}
    g\circ\alpha_H=\alpha'_R\circ \Phi(g)\qquad\forall g\in \Hom_{\cC'}(\beta(H),R).
\end{equation}


 \subsection{$\F_1$-schemes}\vspace{.05in}
 
 We finally recall from \cite{compositio} the notion of $\F_1$-scheme in terms of covariant functors. 
 
 \begin{defn}\label{defnfunfunc} An $\F_1$-functor is a covariant functor from the category
$\Mr=\An\cup_{\beta,\beta^*} \Mo$ to the category of sets.
\end{defn}

The assignment of  an $\F_1$-functor $\mathcal X: \Mr \to \Se$ is equivalent to the specification of the following data:\vspace{.05in}

$\bullet$~An $\Mo$-functor $\underline X$\vspace{.05in}

$\bullet$~A $\Z$-functor $X_\Z$\vspace{.05in}

$\bullet$~A natural transformation $e:\underline X\circ \beta^*\to X_\Z$.\vspace{.05in}

 This natural transformation corresponds to the evaluation of the $\F_1$-functor on the specific morphism  $\alpha'_R$ of \eqref{alphasbis}. 
The notion of $\F_1$-scheme defined in \cite{compositio} is as follows

\begin{defn}\label{defnfonesch} An $\F_1$-scheme is an $\F_1$-functor $\mathcal X: \Mr\longrightarrow \mathfrak{Sets}$, such that:\vspace{.05in}

$\bullet$~The restriction $X_\Z$ of $\mathcal X$ to $\An$ is a $\Z$-scheme.\vspace{.05in}

$\bullet$~The restriction $\underline X$ of $\mathcal X$ to $\Mo$ is an $\Mo$-scheme.\vspace{.05in}

$\bullet$~The natural transformation $e:\underline X\circ \beta^*\to X_\Z$ associated to a field is a bijection (of sets).
\end{defn}

Morphisms of $\F_1$-schemes are natural transformations of the corresponding functors.

\end{document}